\documentclass{amsart}
\usepackage{amscd}

	\title{Multigraded Hilbert Schemes}

	\author{Mark Haiman
	\and    Bernd Sturmfels}

	\address{Department of Mathematics,
	         University of California, 
	         Berkeley, CA 94720}       

	\thanks{Research supported in part by NSF grants DMS-0070772
	(M.H.) and DMS-9970254 (B.S.)}

\theoremstyle{plain}
\newtheorem{theorem}{Theorem}[section]
\newtheorem{lemma}[theorem]{Lemma}
\newtheorem{proposition}[theorem]{Proposition}
\newtheorem{corollary}[theorem]{Corollary}

\theoremstyle{definition}
\newtheorem{definition}[theorem]{Definition}
\newtheorem{remark}[theorem]{Remark}
\newtheorem{example}[theorem]{Example}

\DeclareMathOperator{\Hilb}{Hilb}
\DeclareMathOperator{\Spec}{Spec}
\DeclareMathOperator{\Proj}{Proj}
\DeclareMathOperator{\Hom}{Hom}
\DeclareMathOperator{\Sym}{Sym}
\DeclareMathOperator{\init}{in}
\DeclareMathOperator{\lcm}{lcm}
\DeclareMathOperator{\rank}{rk}

\DeclareMathOperator{\ord}{ord}
\DeclareMathOperator{\red}{red}
\DeclareMathOperator{\integ}{int}
\DeclareMathOperator{\gr}{gr}

\newcommand{\OO}{{\mathcal O}}
\renewcommand{\SS}{{\mathcal S}}

\newcommand{\LL}{{\mathcal L}}

\newcommand{\fun}[1]{\underline{#1}}
\newcommand{\kAlgCateg}{\fun{\text{$k$-Alg}}}

\newcommand{\SetCateg}{\fun{\text{Set}}}
\newcommand{\kSchCateg}{\fun{\text{Sch/$k$}}}
\newcommand{\ZSchCateg}{\fun{\text{Sch/$\Z $}}}
\newcommand{\x}{{\mathbf x}}

\renewcommand{\gg}{{\mathfrak g}}
\newcommand{\slashsub}[1]{/\!_{#1}\,}

\newcommand{\A}{{\mathbb A}}
\renewcommand{\P}{{\mathbb P}}
\newcommand{\C}{{\mathbb C}}
\newcommand{\R}{{\mathbb R}}
\newcommand{\T}{{\mathbb T}}
\newcommand{\Z}{{\mathbb Z}}
\newcommand{\N}{{\mathbb N}}

\begin{document}

\begin{abstract}
We introduce the multigraded Hilbert scheme,
which parametrizes all homogeneous ideals 
with fixed Hilbert function in a polynomial
ring that is graded by any abelian group.
Our construction is widely applicable, it provides
explicit equations, and it allows us to 
prove a range of new results, 
including Bayer's conjecture on equations
defining Grothendieck's classical Hilbert scheme
and the construction of a Chow morphism for
toric Hilbert schemes.

\end{abstract}

\maketitle

\section{Introduction}
\label{sec1}

The multigraded Hilbert scheme parametrizes all ideals in a polynomial
ring which are homogeneous and have a fixed
Hilbert function with respect to a grading by an abelian group.
Special cases include Hilbert schemes of points in affine space
\cite{I}, toric Hilbert schemes \cite{PS}, Hilbert schemes of abelian
group orbits \cite{N}, and Grothendieck's classical Hilbert scheme
\cite{EH}.  We show that the multigraded Hilbert scheme always exists
as a quasiprojective scheme over the ground ring $k$. 
This result is obtained  by means of a
general construction which works in more contexts than just
multigraded polynomial rings.  It also applies to 
$\operatorname{Quot}$ schemes  and to Hilbert schemes arising
in noncommutative geometry; see {\it e.g.},
\cite{AZ}, \cite{BW}.   Our results resolve several open
questions about  Hilbert schemes and their equations.

 Our broader
purpose is to realize the multigraded Hilbert scheme effectively, in
terms of explicit coordinates and defining equations.  These coordinates may
either be global, in the projective case, or local, on affine charts
covering the Hilbert scheme.  
A byproduct of our aim for explicit equations is,
 perhaps surprisingly, a high level of abstract generality.
In particular, we avoid using Noetherian hypotheses, so our results
are valid over any commutative ground ring $k$ whatsoever.

Let $S = k[x_1,\ldots,x_n]$ be the polynomial ring over a commutative
ring $k$.  Monomials $x^u $ in $S$ are identified with vectors $u$ in
$\N^n$.  A grading of $S$ by an abelian group $A$ is a semigroup
homomorphism $\deg \colon \N^n \rightarrow A$. This induces a
decomposition
\[
S = \bigoplus_{a \in A} S_a, \qquad \text{satisfying} \quad S_a \cdot
S_b \subseteq S_{a + b},
\]
where $S_a$ is the $k$-span of all monomials $x^u$ whose degree
is equal to $a$.  
 Note that $S_a$ need not be finitely-generated over $k$.  
 We always assume, without
loss of generality, that the group $A$ is generated by the
elements $a_i = \deg (x_i)$ for $i=1,2, \ldots, n$. 
Let $A_{+} = \deg (\N^{n})$ denote the subsemigroup of $A$ generated by 
$a_1,\ldots,a_n$.

A homogeneous ideal $I$ in $S$ is {\em admissible} if
$(S/I)_a = S_a/I_a$ is a locally free
$k$-module of finite rank (constant on $\Spec k$) for all
$a \in A$.   Its {\em Hilbert function} is 
\begin{equation}\label{e:h-of-I}
h_I \colon A \rightarrow \N, \quad h_{I}(a) = \rank _{k} (S/I)_a.
\end{equation}
Note that the support of $h_{I}$ is necessarily contained in $A_{+}$.
Fix any numerical function $h \colon A \rightarrow \N $ supported on
$A_{+}$.  We shall construct a scheme over $k$ which parametrizes, in
the technical sense below, all admissible ideals $I$ in $S$ with $h_I
= h$.

Recall ({\it e.g.} from \cite{EH}) that every scheme $Z$ over $k$ is
characterized by its {\em functor of points}, which maps the category
of $k$-algebras to the category of sets as follows:
\begin{equation}\label{e:scheme-functor}
\fun{Z}\colon \kAlgCateg \rightarrow \SetCateg, \quad \fun{Z}(R) =
\Hom(\Spec R, Z).
\end{equation}
Given our graded polynomial ring $S = k[x_1,\ldots,x_n]$ and
Hilbert function $h$, the {\em Hilbert functor} $H_{S}^{h}\colon
\kAlgCateg \rightarrow \SetCateg $ is defined as follows:
$H_{S}^{h}(R)$ is the set of homogeneous 
ideals $I\subseteq {R\otimes_k S}$ such
that $({R\otimes S_{a}})/I_{a}$ is a locally free $R$-module of rank
$h(a)$ for each $a\in A$.  We shall construct the scheme which
represents this functor.

\begin{theorem}\label{t:multi-Hilb}
There exists a quasiprojective scheme $Z$ over $k$ such that $\fun{Z}
= H_{S}^{h}$.
\end{theorem}

The scheme $Z$ is called the {\em multigraded Hilbert scheme} and is
also denoted $ H_{S}^{h}$.   It is projective if the grading is {\em
positive}, which means that $x^{0}=1$ is the only monomial of degree
$0$.  Note that if the grading is positive, then $A_{+}\cap -(A_{+}) =
\{0 \}$.

\begin{corollary}
\label{p:projective} If the grading of the polynomial ring $S =
k[x_1,\ldots,x_n]$ is positive then the multigraded Hilbert scheme
$H_{S}^{h}$ is projective over the ground ring $k$.
\end{corollary}

This corollary also follows from recent work of Artin and Zhang
\cite{AZ}.  The approach of Artin and Zhang is non-constructive, and
does not apply when the $S_{a}$ are not finite over $k$ and the
Hilbert scheme is only quasiprojective.  In our approach, the
Noetherian and finite-generation hypotheses in \cite{AZ} are replaced
by more combinatorial conditions.  This gives us sufficient generality
to construct quasiprojective Hilbert schemes, and the proof becomes
algorithmic, transparent and uniform, requiring no restrictions on the
ground ring $k$, which need not even be Noetherian.

This paper is organized as follows.  In Section~\ref{sec2} we present
a general construction realizing Hilbert schemes as quasiprojective
varieties. The main results in Section~\ref{sec2} are
Theorems~\ref{t:Hilb-scheme-I} and \ref{t:Hilb-scheme-II}.  In Section
\ref{sec3}, we apply these general theorems to prove Theorem
\ref{t:multi-Hilb} and Corollary \ref{p:projective}.  The needed
finiteness hypotheses are verified using Maclagan's finiteness theorem
\cite{M} for monomial ideals in $S$.  Our main results in Section
\ref{sec3} are Theorems \ref{t:multi-iso} and \ref{c:determinants}.
These two theorems identify finite subsets $D$ of the group $A$ such
that the degree restriction morphism $H^h_{S} \rightarrow H^h_{S_D}$
is a closed embedding (respectively an isomorphism), and they lead to
explicit determinantal equations and quadratic equations for the
Hilbert scheme $H_{S}^{h}$.

Section~\ref{sec4} concerns the classical Grothendieck Hilbert scheme
which parametrizes ideals with a given Hilbert polynomial 
(as opposed to a given Hilbert function) in the usual $\N$-grading.  
The results of Gotzmann \cite{Got} can be
interpreted as identifying the Grothendieck Hilbert scheme
with our $H_{S}^{h}$, for a suitably chosen Hilbert function
$h$, depending on the Hilbert polynomial.  Our construction 
naturally yields two descriptions of the Hilbert scheme by coordinates
and equations.  The first reproduces Gotzmann's equations in terms of
Pl\"ucker coordinates in two consecutive degrees.  The second
reproduces equations in terms of Pl\"ucker coordinates in just one
degree. We prove a conjecture from Bayer's 1982 thesis \cite{Bay}
stating that Bayer's set-theoretic equations of degree $n$
actually define the Hilbert scheme as a scheme.

In Section~\ref{sec5} we examine the case where $h$ is the incidence
function of the semigroup  $A_+$, in which case 
$H_{S}^{h}$ is called the {\em toric Hilbert
scheme}.  In the special cases when the grading is positive or when
the group $A$ is finite, this scheme was constructed by Peeva and
Stillman \cite{PS} and Nakamura \cite{N} respectively.  We unify and
extend results by these authors, and we resolve Problem 6.4 in \cite{S}
by constructing the natural morphism from the toric Hilbert scheme to
the {\em toric Chow variety}.

Recent work by Santos \cite{San} provides an example where both the
toric Chow variety and the toric Hilbert scheme are {\em
disconnected}. This shows that the multigraded Hilbert scheme
$H_{S}^{h}$ can be disconnected, in contrast to Hartshorne's classical
connectedness result \cite{Har} for the Grothendieck Hilbert scheme.

In Section~\ref{sec6} we demonstrate that the results of
Section~\ref{sec2} are applicable to a wide range of parameter spaces
other than the multigraded Hilbert scheme; specifically, we construct
$\operatorname{Quot}$ schemes and Hilbert schemes parametrizing ideals
in the Weyl algebra, the exterior algebra and other noncommutative rings.  

\smallskip 

Before diving into the abstract setting of Section \ref{sec2},
we wish to first present a few concrete examples and basic facts
concerning multigraded Hilbert schemes.

\smallskip 

\begin{example} \label{Ex1}
Let $n = 2 $ and $k = \C$, the complex numbers, and fix $S = \C[x,y]$.
We conjecture that $H_{S}^{h}$ is smooth and irreducible for any 
group $A$ and any $h : A \rightarrow \N$.
\begin{itemize}

\item[(a)] If $A = 0$ then $H_{S}^{h}$ is the Hilbert scheme 
of $n = h(0)$ points in the affine plane $\A^2$.  This scheme
is smooth and irreducible of dimension $2n$; see \cite{Fog}.

\item[(b)] If $A=\Z$, $\deg(x)$, $\deg(y)$ are positive integers, and $h$
has finite support, then $H_{S}^{h}$ is an irreducible
component in the  fixed-point set of a $\C^{*}$-action on 
the Hilbert scheme of points; see e.g.~\cite{E}.
This was proved by Evain \cite{E2}.

\item[(c)] If $A = \Z$, $\deg(x) = \deg(y) = 1$ and $h(a) = 1$ for $a \geq
0$, then   $H_{S}^{h} = \P^1$.

\item[(d)] More generally, if $A = \Z$, $\deg(x) = \deg(y) = 1$ and
$h(a)=\min (m,a+1)$, for some integer $m \geq 1$, 
 then $H_{S}^{h}$ is the Hilbert scheme of $m$
points on $\P ^{1}$.

\item[(e)] If $A = \Z$, $\deg(x) = - \deg(y) = 1$ and $h(a) = 1$ for all
$a$, then $H_{S}^{h} = \A^1$.

\item[(f)] If $A = \Z^2$, $\deg(x) = (1,0)$ and $\deg(y) = (0,1)$, then
$H_{S}^{h}$ is either empty or a point. In the latter
case it consists of a single monomial ideal.

\item[(g)] If $A = \Z/2\Z$, $\deg(x) = \deg(y) = 1$ and $h(0) = h(1) = 1$,
then $H_{S}^{h}$ is isomorphic to the cotangent bundle of the projective line
$\P^{1}$.
\end{itemize}
\end{example}

\smallskip 

\begin{example} \label{Ex2}
Let $n = 3$. This example is  the  smallest
reducible Hilbert scheme known to us. 
We fix the $\Z^2$-grading of the polynomial ring $S = \C[x,y,z]$ given by  
$$  {\rm deg}(x) = (1,0) \,, \,\,\, {\rm deg}(y) = (1,1)
 \,, \,\,\, {\rm deg}(z) = (0,1) . $$
Consider the closed subscheme in the Hilbert scheme of nine points 
in  $\A^3$ consisting of homogeneous ideals $I \subset S$
such that $S/I$ has the bivariate Hilbert series $$ s^2 t^2 \, +
 \, s^2 t \, + \, s t^2 \, + \, s^2 \, +2 s t \, + \, s \, +\,t \, +\, 1 . $$
This multigraded Hilbert scheme  is the reduced union of two projective lines 
$\P^{1}$
which intersect in a common torus fixed point. The universal family equals
$$ \langle \,  x^3,  \,  x y^2,  \,  x^2 y, \,
 y^3, \,  a_0 x^2 z - a_1 x y,  \,  b_0 x y z - b_1 y^2 , \,  y^2 z,  \,
 z^2  \, \rangle \quad \hbox{with} \quad \,  a_1 b_1 \,\, = \,\, 0. $$
Here $(a_0:a_1)$ and $(b_0:b_1)$ are coordinates on  two projective lines.
This Hilbert scheme has three torus fixed points, namely,
the three monomial ideals in the family.
\end{example}

\smallskip

In these examples we saw that if 
the Hilbert function $h$ has finite support, say
$\,m = \sum_{a \in A} h(a) $, then $\,H_S^h \,$
is a closed subscheme of the
Hilbert scheme of $m$ points in $\A^n$.
More generally,
there is a canonical embedding of one multigraded Hilbert scheme into
another when the grading and Hilbert function of the first refine
those of the second.  Let $\phi \colon A_{0}\rightarrow A_{1}$ be a
homomorphism of abelian groups.  A grading $\deg _{0}\colon \N
^{n}\rightarrow A_{0}$ {\em refines} $\deg _{1}\colon \N
^{n}\rightarrow A_{1}$ if $\deg_{1} = \phi \circ \deg _{0}$.  In this
situation, a function $h_{0}\colon A_{0}\rightarrow \N $
{\em refines } $h_{1}\colon A_{1}\rightarrow \N $ if $h_{1}(u) = \sum
_{\phi (v) =u} h_{0}(v)$ for all $u\in A_{1}$. Any admissible ideal 
$I\subseteq {R\otimes S}$ with Hilbert function $h_{0}$ for the grading 
$\deg _{0}$ is also admissible with Hilbert function $h_{1}$ for $\deg
_{1}$.  Hence the Hilbert functor $H^{h_{0}}_{S}$ is a subfunctor
of $H^{h_{1}}_{S}$. The following assertion
will be proved in Section \ref{sec3}.

\smallskip

\begin{proposition}\label{p:refinement}
If $(\deg _{0}, h_{0})$ refines $(\deg _{1}, h_{1})$, then the natural
embedding of Hilbert functors described above is induced by an
embedding of the multigraded Hilbert scheme 
$H_{S}^{h_0}$ as a closed subscheme of $H_{S}^{h_1}$.
\end{proposition}

\smallskip

A nice feature of the multigraded Hilbert scheme, in common with other
Hilbert schemes, is that its tangent space at any point has a simple
description.  We assume that $k$ is a field and $I \in H_S^h(k)$.  The
$S$-module ${\rm Hom}_S (I, S/I)$ is graded by the group $A$, and each
component $({\rm Hom}_S (I, S/I))_a$ is a finite-dimensional
$k$-vector space.

\begin{proposition}
For $k$ a field, the Zariski tangent space to the multigraded Hilbert
scheme $H_{S}^h$ at a point $I\in H_{S}^{h}(k)$ is  canonically
isomorphic to $\,({\rm Hom}_S (I, S/I))_0$.
\end{proposition}

\begin{proof}
Let $R = k[\epsilon]/\langle \epsilon^2 \rangle $. The tangent space
at $I$ is the set of points in $H_S^h(R)$ whose image under the
map $H_S^h(R) \rightarrow H_S^h(k)$ is $I$. Such a point is an
$A$-homogeneous ideal $J \subset R[\x ]= k[\x,\epsilon]/\langle
\epsilon^2 \rangle$ such that $J / \langle \epsilon \rangle $ equals
the ideal $I$ in $S = k[\x]$ and $R[\x ]/J$ is a free $R$-module.
Consider the map from $ k[\x]$ to $ \epsilon R[\x] \cong  k[\x]$ given
by multiplication by $\epsilon$.  This multiplication map followed by
projection onto $\epsilon R[\x]/ (J\cap \epsilon R[\x]) \cong 
k[\x]/I$ represents a degree zero homomorphism $ I \rightarrow S/I$,
and, conversely, every degree zero homomorphism $ I \rightarrow S/I$
arises in this manner from some $J$.
\end{proof}

\section{A general framework for  Hilbert schemes}
\label{sec2}

Fix a commutative ring $k$ and an arbitrary index set $A$ called
``degrees.''  Let
\begin{equation}\label{e:T}
T \; = \; \bigoplus _{a\in A} T_{a}
\end{equation}
be a graded $k$-module, equipped with a collection of operators $F =
\bigcup _{a,b\in A} F_{a,b}$, where $F_{a,b}\subseteq \Hom
_{k}(T_{a},T_{b})$.  Given a commutative $k$-algebra $R$, we denote by
${R\otimes T}$ the graded $R$-module $\bigoplus _{a} {R\otimes
T_{a}}$, with operators $\hat{F}_{a,b} = ({1_{R}\otimes -})(F_{a,b})$.
A homogeneous submodule $L = \bigoplus _{a} L_{a} \subseteq {R\otimes
T}$ is an {\em $F$-submodule} if it satisfies
$\hat{F}_{a,b}(L_{a})\subseteq L_{b}$ for all $a,b\in A$.  We may
assume that $F$ is closed under composition: $F_{bc}\circ
F_{ab}\subseteq F_{ac}$ for all $a,b,c\in A$ and $F_{aa}$ contains the
identity map on $T_{a}$ for all $a\in A$.  In other words, $(T,F)$ is
a small category of $k$-modules, with the components $T_{a}$ of $T$ as
objects and the elements of $F$ as arrows.

Fix a function $h\colon A\rightarrow \N $. Let $H_{T}^{h}(R)$ be the
set of $F$-submodules $L\subseteq R\otimes T$ such that $(R\otimes
T_{a})/L_{a}$ is a locally free $R$-module of rank $h(a)$ for each
$a\in A$.  If $\phi \colon R\rightarrow S$ is a homomorphism of
commutative rings (with unit), then local freeness implies that $L'=
{S\otimes _{R}L}$ is an $F$-submodule of ${S\otimes T}$, and $({S
\otimes T_{a}})/L'_{a}$ is locally free of rank $h(a)$ for each $a\in
A$. Defining $H_{T}^{h}(\phi )\colon H_{T}^{h}(R) \rightarrow
H_{T}^{h}(S)$ to be the map sending $L$ to $L'$ makes $H_{T}^{h}$ a
functor $ \kAlgCateg \rightarrow \SetCateg $, called the {\em Hilbert
functor}.

If $(T,F)$ is a graded $k$-module with operators, as above, and
$D\subseteq A$ is a subset of the degrees, we denote by
$(T_{D},F_{D})$ the restriction of $(T,F)$ to degrees in $D$.  In the
language of categories, $(T_{D},F_{D})$ is the full subcategory of
$(T,F)$ with objects $T_{a}$ for $a\in D$.  There is an obvious
natural transformation of Hilbert functors $H_{T}^{h}\rightarrow
H_{T_{D}}^{h}$ given by restriction of degrees, that is, $L\in
H_{T}^{h}(R)$ goes to $L_{D} = \bigoplus _{a\in D}L_{a}$.

\begin{remark}\label{r:Fclosed} Given an $F_{D}$-submodule $L\subseteq
{R\otimes T_{D}}$, let $L'\subseteq {R\otimes T}$ be the $F$-submodule
it generates.  The assumption that $F$ is closed under composition
implies that $L'_{a} = \sum _{b\in D} F_{ba}(L_{b})$.  In particular,
the restriction $L'_{D}$ of $L'$ is equal to $L$.
\end{remark}

We show that, under suitable hypotheses, the Hilbert functor
$H_{T}^{h}$ is represented by a quasiprojective scheme over $k$,
called the {\em Hilbert scheme}.  Here and elsewhere we will abuse
notation by denoting this scheme and the functor it represents by the
same symbol, so we also write $H_{T}^{h}$ for the Hilbert scheme.

\begin{theorem}\label{t:Hilb-scheme-I}
Let $(T,F)$ be a graded $k$-module with operators, as above.  Suppose   
that $M\subseteq N\subseteq T$ are homogeneous $k$-submodules satisfying     
four conditions:
\begin{itemize}                                                         
\item [(i)] $N$ is a finitely generated $k$-module;                     
\item [(ii)] $N$ generates $T$ as an $F$-module;
\item [(iii)] for every field $K\in \kAlgCateg $ and every               
$L\in H_{T}^{h}(K)$, $M$ spans $({K\otimes T}) / L$;  and
\item [(iv)] there is a subset $G\subseteq F$, generating $F$ as a      
category, such that $GM\subseteq N$.                                
\end{itemize}                                                           
Then $H_{T}^{h}$ is represented by a quasiprojective scheme over $k$.
It is a closed subscheme of the relative Grassmann scheme
$G^{h}_{N\setminus M}$, which is defined below.
\end{theorem}

In hypothesis (iii), $N$ also spans $({K\otimes T})/L$, so $\dim_{K}
({K\otimes T})/L = \sum _{a\in A}h(a) $ is finite.  Therefore
Theorem~\ref{t:Hilb-scheme-I} only applies when $h$ has finite
support. Our strategy in the general case is to construct the Hilbert
scheme for a finite subset $D$ of the degrees $A$ and then to use the
next theorem to refine it to all degrees.

\begin{theorem}\label{t:Hilb-scheme-II}                          
Let $(T,F)$ be a graded $k$-module with operators and 
$D \subseteq A$ such that $H_{T_{D}}^{h}$ is represented by a scheme
over $k$.  Assume that for each degree $a\in A$:
\begin{itemize}                                                       
\item [(v)] there is a finite subset $E\subseteq \bigcup _{b\in D}
F_{ba}$ such that $T_{a}/\sum _{b\in D} E_{ba}(T_{b})$ is a
finitely generated $k$-module; and
\item [(vi)] for every field $K\in \kAlgCateg $ and every $L_{D}\in
H_{T_{D}}^{h}(K)$, if $L'$ denotes the $F$-submodule of ${K\otimes T}$
generated by $L_{D}$, then $\dim ({K\otimes T_{a}})/L'_{a}\leq h(a)$.
\end{itemize}                                                         
Then the natural transformation $H_{T}^{h}\rightarrow H_{T_{D}}^{h}$  
makes $H_{T}^{h}$ a subfunctor of $H_{T_{D}}^{h}$, represented by a   
closed subscheme of the Hilbert scheme $H_{T_{D}}^{h}$.               
\end{theorem}                                                         

We realize that conditions (i)--(vi) above appear obscure at first
sight.  Their usefulness will become clear as we apply these theorems
in Section~\ref{sec3}.

Sometimes the Hilbert scheme is not only quasiprojective over $k$, but
projective.

\begin{corollary} \label{SometimesProjective}
In Theorem \ref{t:Hilb-scheme-I}, in place of hypotheses {\rm
(i)--(iv)}, assume only that the degree set $A$ is finite, and $T_a$
is a finitely-generated $k$-module for all $a \in A$.  Then
$H_{T}^{h}$ is projective over $k$.
\end{corollary}

\begin{proof}
We can take $M = N = T$ and $G=F$.  Then hypotheses (i)--(iv) are
trivially satisfied, and the relative Grassmann scheme
$G^{h}_{N\setminus M}$ in the conclusion is just the Grassmann scheme
$G^{h}_{N}$.  It is projective by Proposition
\ref{p:GNr-Grass-scheme}, below.
\end{proof}

\begin{remark}\label{r:projective}
In Theorem~\ref{t:Hilb-scheme-II}, suppose in addition to hypotheses
{\rm (v)} and {\rm (vi)} that $D$ is finite and $T_{a}$ is finitely
generated for all $a\in D$.  Then we can again conclude that
$H_{T}^{h}$ is projective, since it is a closed subscheme of the
projective scheme $H_{T_{D}}^{h}$.
\end{remark}

In what follows we review some facts about functors, Grassmann
schemes, and the like, then turn to the proofs of
Theorems~\ref{t:Hilb-scheme-I} and \ref{t:Hilb-scheme-II}.  In Section
3 we use these theorems to construct the multigraded Hilbert scheme.

We always work in the category $\kSchCateg $ of schemes over a fixed
ground ring $k$.  We denote the functor of points of a scheme $Z$ by
$\fun{Z}$ as in \eqref{e:scheme-functor}.

\begin{proposition}[{\cite[Proposition VI-2]{EH}}]
\label{p:functor-determines-scheme} The scheme $Z$ is characterized by
its functor $\fun{Z}$, in the sense that every natural transformation
of functors $\fun{Y}\rightarrow \fun{Z}$ is induced by a unique
morphism $Y\rightarrow Z$ of schemes over $k$.
\end{proposition}

Our approach to the construction of Hilbert schemes will be to
represent the functors in question by subschemes of Grassmann schemes.
The theoretical tool we need for this is a representability theorem
for a functor defined relative to a given scheme functor. The
statement below involves the concepts of {\em open subfunctor}, see
\cite[\S VI.1.1]{EH}, and {\em Zariski sheaf}, introduced as ``sheaf
in the Zariski topology'' at the beginning of \cite[\S
VI.2.1]{EH}. Being a Zariski sheaf is a necessary condition for a
functor $\kAlgCateg \rightarrow \SetCateg $ to be represented by a
scheme. See \cite[Theorem VI-14]{EH} for one possible converse. Here
is the relative representability theorem we will use.

\begin{proposition}\label{p:relative-representability}
Let $\eta \colon Q\rightarrow \fun{Z}$ be a natural transformation of
functors $\kAlgCateg \rightarrow \SetCateg $, where $\fun{Z}$ is a
scheme functor and $Q$ is a Zariski sheaf.  Suppose that $Z$ has a
covering by open sets $U_{\alpha }$ such that each subfunctor $\eta
^{-1}(\fun{U_{\alpha }})\subseteq Q$ is a scheme functor.  Then $Q$ is
a scheme functor, and $\eta $ corresponds to a morphism of
schemes.
\end{proposition}

\begin{proof}
Let $Y_{\alpha }$ be the scheme whose functor is $\eta
^{-1}(\fun{U_{\alpha }})$.  The induced natural transformation $\eta
^{-1}(\fun{U_{\alpha }})\rightarrow \fun{U_{\alpha }}$ provides us
with a morphism $\pi _{\alpha }\colon Y_{\alpha }\rightarrow U_{\alpha
}$.  For each $\alpha $ and $\beta $, the open subscheme $\pi _{\alpha
}^{-1}(U_{\alpha }\cap U_{\beta })\subseteq Y_{\alpha }$ has functor
$\eta ^{-1}(\fun{U_{\alpha }}\cap \fun{U_{\beta }})$.  In particular,
we have a canonical identification of $\pi _{\alpha }^{-1}(U_{\alpha
}\cap U_{\beta })$ with $\pi _{\beta }^{-1}(U_{\alpha }\cap U_{\beta
})$, and these identifications are compatible on every triple
intersection $U_{\alpha }\cap U_{\beta }\cap U_{\gamma}$.  By the
gluing lemma for schemes, there is a scheme $Y$ with a morphism $\pi
\colon Y\rightarrow Z$ such that for each $\alpha $ we have $Y_{\alpha
} = \pi ^{-1}(U_{\alpha })$ and $\pi _{\alpha } = \pi |_{Y_{\alpha
}}$.

Let $R$ be a $k$-algebra and let $\phi $ be an element of
$\fun{Y}(R)$, that is, a morphism $\phi \colon \Spec R\rightarrow Y$.
Since the $Y_{\alpha }$ form an open covering of $Y$, there are
elements $f_{i}$ generating the unit ideal in $R$ such that $\phi $
maps $U_{f_{i}}\subseteq \Spec R$ into some $Y_{\alpha _{i}}$.  Let
$\phi _{i}\colon U_{f_{i}}\rightarrow Y_{\alpha _{i}}$ be the
restriction of $\phi $; it is an element of $\fun{Y_{\alpha
_{i}}}(R_{f_{i }})\subseteq Q(R_{f_{i}})$.  For each $i$, $j$, the
elements $\phi _{i}$, $\phi _{j}$ restrict to the same morphism $\phi
_{ij}\colon U_{f_{i}f_{j}}\rightarrow Y_{\alpha _{i}}\cap Y_{\alpha
_{j}}$, and therefore have the same image in $Q(R_{f_{i}f_{j}})$.
Since $Q$ is a Zariski sheaf by hypothesis, the elements $\phi _{i}$
are all induced by a unique element $\hat{\phi }\in Q(R)$.

We have thus constructed a transformation $\xi \colon
\fun{Y}\rightarrow Q$ sending $\phi \in \fun{Y}(R)$ to $\hat{\phi }\in
Q(R)$, and it is clearly natural in $R$.  We claim that $\xi $ is a
natural isomorphism.  First note that $\hat{\phi }$ determines each
$\phi _{i}$ by definition, and these determine $\phi $ since the
$U_{f_{i}}$ cover $\Spec R$.  Hence $\xi _{R}$ is injective.  Now
consider any $k$-algebra $R$ and $\lambda \in Q(R)$.  Then $\eta
(\lambda )\in \fun{Z}(R)$ is a morphism $ \Spec R\rightarrow Z$, and
we can cover $\Spec R$ by open sets $U_{f_{i}}$ such that $\eta
(\lambda )$ maps each $U_{f_{i}}$ into some $U_{\alpha _{i}}$.  This
means that the image of $\lambda $ in $Q(R_{f_{i}})$ belongs to $\eta
^{-1}(\fun{U_{\alpha _{i}}})$, that is, to $\fun{Y_{\alpha _{i}}}$.
Since $\fun{Y}$ is a Zariski sheaf and the $U_{f_{i}}$ cover $\Spec
R$, this implies that $\lambda $ belongs to $\xi _{R} (\fun{Y}(R))$.
Hence $\xi $ is surjective.
\end{proof}

\begin{corollary}\label{c:closed-subfunctor}
Under the hypotheses of Proposition \ref{p:relative-representability},
if the natural transformations $\eta ^{-1}(\fun{U_{\alpha
}})\rightarrow \fun{U_{\alpha }}$ given by restricting $\eta $ are
induced by closed embeddings of schemes, then so is $\eta $.
\end{corollary}

\begin{proof}
This just says that the condition for a morphism $\eta \colon
Y\rightarrow Z$ to be a closed embedding is local on $Z$.  Indeed, the
result is valid with ``closed embedding'' replaced by any property of
a morphism that is local on the base.
\end{proof}

Another useful characterization of natural transformations $\eta
\colon Q\rightarrow \fun{Z}$ represented by closed subschemes of $Z$
is as subfunctors defined by a closed condition.  A condition on
$R$-algebras is {\em closed} if there exists an ideal $I\subseteq R$
such that the condition holds for an $R$-algebra $S$ if and only if
the image of $I$ in $S$ is zero.

Let $Z$ be a scheme over $k$ and $\eta \colon Q\hookrightarrow
\fun{Z}$ a subfunctor.  We wish to decide whether $\eta $ is
represented by a closed embedding.  Consider a $k$-algebra $R$ and an
element $\lambda \in \fun{Z}(R)$, or equivalently a morphism $\lambda
\colon \Spec R\rightarrow Z$.  Given this data, we define a condition
$V_{R,\lambda }$ on $R$-algebras $S$, as follows.  Let $\phi \colon
R\rightarrow S$ be the ring homomorphism making $S$ an $R$-algebra.
Then $S$ satisfies the condition $V_{R,\lambda }$ if the element
$\fun{Z}(\phi )\lambda \in \fun{Z}(S)$ belongs to the subset $\eta
_{S}( Q(S))\subseteq \fun{Z}(S)$ defined by the subfunctor.  We can
now express the content of Proposition
\ref{p:relative-representability} and Corollary
\ref{c:closed-subfunctor} as follows.

\begin{proposition}\label{p:closed-condition}
Let $\eta \colon Q\hookrightarrow \fun{Z}$ be a subfunctor, where
$\fun{Z}$ is a scheme functor and $Q$ is a Zariski sheaf.  Then $Q$ is
represented by a closed subscheme of $Z$ if and only if $V_{R,\lambda
}$ is a closed condition for all $R\in \kAlgCateg $ and $\lambda \in
\fun{Z}(R)$.
\end{proposition}

\begin{proof}
First suppose that $Y\subseteq Z$ is a closed subscheme, and $Q=
\fun{Y}$ is the corresponding subfunctor of $\fun{Z}$.  Given $\lambda
\colon \Spec R\rightarrow Z$, let $I\subseteq R$ be the ideal defining
the scheme-theoretic preimage $\lambda ^{-1}(Y)\subseteq \Spec R$.
The condition $V_{R,\lambda }$ on an $R$-algebra $S$ is that $\phi
\colon R\rightarrow S$ factor through $R/I$, so it is a closed
condition.

For the converse, using Proposition \ref{p:relative-representability}
and Corollary \ref{c:closed-subfunctor}, it suffices to verify that
$Q' = Q\cap \fun{U}$ is represented by a closed subscheme of $U$, for
each $U = \Spec R$ in an affine open covering of $Z$.  The inclusion
$\lambda \colon U\hookrightarrow Z$ is an element $\lambda \in
\fun{Z}(R)$.  The subset $Q'(S)\subseteq \fun{U}(S)$ is the set of
morphisms $\nu \colon \Spec S\rightarrow U$ such that $\lambda \circ
\nu $ belongs to $\eta_{S} (Q(S))$.  If $\phi \colon R\rightarrow S$
is the ring homomorphism underlying such a morphism $\nu $, then
$\lambda \circ \nu = \fun{Z}(\phi )\lambda$, so $\nu $ belongs to
$Q'(S)$ if and only if the $R$-algebra $S$ satisfies the condition
$V_{R,\lambda }$. By hypothesis, the closed condition $V_{R,\lambda }$
is defined by an ideal $I\subseteq R$. Hence $Q'(S)$ is naturally
identified with the set of ring homomorphisms $\phi \colon
R\rightarrow S$ that factor through $R/I$.  In other words, $Q'$ is
represented by the closed subscheme $V(I)\subseteq U = \Spec R$.
\end{proof}

Recall that an $R$-module $W$ is {\em locally free of rank $r$} if
there exist $f_{1},\ldots,f_{k}\in R$ generating the unit ideal, such
that $W_{f_{i}}\cong R_{f_{i}}^{r}$ for each $i$.  Let $N$ be any
finitely generated $k$-module.  The Grassmann scheme $G_{N}^{r}$
represents the {\em Grassmann functor}, defined as follows: for $R\in
\kAlgCateg $, the set $G_{N}^{r}(R)$ consists of all submodules
$L\subseteq {R\otimes N}$ such that $({R\otimes N})/L$ is locally free
of rank $r$.

We review the description of the Grassmann scheme $G_{N}^{r}$ in terms
of coordinates, starting with the free module $N = k^m$, whose basis
we denote by $X$.  For this $N$ we write $G_{m}^{r}$ in place of
$G_{N}^{r}$.  Consider a subset $B\subseteq X$ with $r$ elements.  Let
$G_{m\setminus B}^{r}\subseteq G_{m}^{r}$ be the subfunctor describing
submodules $L\in R^{m}$ such that $R^{m}/L$ is free with basis $B$.
This subfunctor is represented by the affine space $\A ^{r(m-r)} =
\Spec k[\gamma ^{x}_{b}:x\in X \backslash B, b\in B]$.  Evaluated at
$L\in G_{m\setminus B}^{r}(R)$, the coordinate $\gamma^{x}_{b}\in R$
is given by the coefficient of the basis vector $b$ in the unique
expansion of $x$ modulo $L$.  We also set $\gamma^{x}_{b} = \delta
_{x,b}$ for $x\in B$.  Passing to Pl\"ucker coordinates, one proves
(see \cite[Exercise VI-18]{EH}) that the Grassmann functor $G_{m}^{r}$
is represented by a projective scheme over $k$, called the {\em
Grassmann scheme}, and the subfunctors $G_{m\setminus B}^{r}$ are
represented by open affine subsets which cover the Grassmann scheme
$G_m^r$.

Next consider an arbitrary finitely-generated $k$-module $N =
k^{m}/J$.  For any $k$-algebra $R$, the module ${R\otimes N}$ is 
isomorphic to $R^{m}/RJ$.  The Grassmann
functor $G^{r}_{N}$ is naturally isomorphic to the subfunctor of
$G^{r}_{m}$ describing submodules $L\subseteq R^{m}$ such that
$RJ\subseteq L$.  If $R^{m}/L$ has basis $B \subseteq X$, then the
condition $RJ\subseteq L$ can be expressed as follows: for each $u\in
J$, write $u = \sum _{x\in X}{a^{u}_{x} \cdot x}$, with $a^{u}_{x}\in
k$.  Then
\begin{equation}\label{e:J-in-L-condition}
\sum _{x\in X} a^{u}_{x} \cdot \gamma ^{x}_{b} = 0 \quad \text{for all
$u \in J$ and $b\in B$.}
\end{equation}
It follows that, for each $B$, the intersection of subfunctors
$G^{r}_{m\setminus B}\cap G^{r}_{N}\subseteq G^{r}_{m}$ is represented
by the closed subscheme of $\Spec k[\gamma ^{x}_{b}]$ defined by the
$k$-linear equations in \eqref{e:J-in-L-condition}.  The condition
$RJ\subseteq L$ is local on $R$, so the subfunctor $G_{N}^{r}\subseteq
G_{m}^{r}$ is a Zariski sheaf.  Therefore Proposition
\ref{p:relative-representability} and Corollary
\ref{c:closed-subfunctor} give the following result.

\begin{proposition}\label{p:GNr-Grass-scheme}
Let $N$ be a finitely generated $k$-module.  The Grassmann functor
$G_{N}^{r}$ is represented by a closed subscheme of the classical
Grassmann scheme $G_{m}^{r}$, called the {\em Grassmann scheme of
$N$}.  In particular, it is projective over $k$.
\end{proposition}

Now suppose that we are given a submodule $M\subseteq N$ (not
necessarily finitely generated, as we are not assuming $k$ is
Noetherian).  For any set $B$ of $r$ elements in $M$, we can choose a
presentation of $N$ in which the generators $X$ contain $B$. The
intersection of $G^{r}_{N}$ with the standard open affine
$G^{r}_{m\setminus B}$ defines an open affine subscheme
$G^{r}_{N\setminus B}\subseteq G^{r}_{N}$.  The affine scheme
$G^{r}_{N\setminus B}$ parametrizes quotients $({R\otimes N})/L$ that
are free with basis $B$.  The union of the subschemes
$G^{r}_{N\setminus B}$ over all $r$-element subsets $B\subseteq M$ is
an open subscheme $G^{r}_{N\setminus M}$ of $G^{r}_{N}$.  The
corresponding subscheme functor describes quotients $({R\otimes N})/L$
that are locally free with basis contained in $M$.  In other words,
$L\in G^{r}_{N}(R)$ belongs to $G^{r}_{N\setminus M}(R)$ if and only
if there are elements $f_{1},\ldots,f_{k}$ generating the unit ideal
in $R$, such that each $({R\otimes N/L})_{f_{i}}$ has a basis
$B_{i}\subseteq M$.  Equivalently, $L$ belongs to $G^{r}_{N\setminus
M}(R)$ if and only if $M$ generates $({R\otimes N})/L$, since the
latter is a local condition on $R$.  The subfunctor $G^{r}_{N\setminus
M}$ of the Grassmann functor $G^{r}_{N}$ is called the {\em relative
Grassmann functor}.

\begin{proposition}\label{p:relative-Grassmann}
Let $N$ be a finitely generated $k$-module and $M$ a submodule.  The
functor $G^{r}_{N\setminus M}$ is represented by an open subscheme of
$G^{r}_{N}$, called the {\em relative Grassmann scheme of $M \subseteq
N$}.  In particular, it is quasiprojective over $k$.
\end{proposition}

Note that if $M = N$ then the relative Grassmann scheme
$G^{r}_{N\setminus M}$ coincides with $G^{r}_N$ and is therefore
projective.  If $M$ is any submodule of $N$ then the open subscheme
$G^{r}_{N\setminus M}\subseteq G^{r}_{N}$ can be described in local
affine coordinates as follows.  Fix a set of $r$ elements $B\subseteq
N$ and consider the standard affine in $G^{r}_{N}$ describing
submodules $L$ such that $({R\otimes N})/L$ has basis $B$.  We form a
matrix $\Gamma $ with $r$ rows, and columns indexed by elements $x\in
M$, whose entries in each column are the coordinate functions $\gamma
^{x}_{b}$ for $b\in B$.  Then $G^{r}_{N\setminus M}$ is described
locally as the complement of the closed locus defined by the vanishing
of the $r\times r$ minors of $\Gamma $.

The definitions and results on Grassmann schemes extend readily to
homogeneous submodules of a finitely generated graded module $N =
\bigoplus _{a\in A} N_{a}$, where $A$ is a finite set of ``degrees.''
Fix a function $h\colon A\rightarrow \N $.  We define the graded
Grassmann functor $G_{N}^{h}$ by the rule that $G_{N}^{h}(R)$ is the
set of homogeneous submodules $L\subseteq {R\otimes N}$ such that
$({R\otimes N_{a}})/L_{a}$ is locally free of rank $h(a)$ for all
$a\in A$.  To give such a submodule $L$, it is equivalent to give each
$L_{a}$ separately.  Thus $G_{N}^{h}$ is naturally isomorphic to the
product $\prod _{a\in A}G_{N_{a}}^{h(a)}$, and in particular it is
projective over $k$.  Similarly, the {\em relative graded Grassmann
functor} $G_{N\setminus M}^{h}$, where $M\subseteq N$ is a homogeneous
submodule, is represented by a quasiprojective scheme over $k$.

\begin{remark}\label{r:refinement}
In the graded situation, $G^{h}_{N}$ is a subfunctor of the ungraded
Grassmann functor $G^{r}_{N}$, where $r = \sum _{a}h(a)$.  Similarly,
$G_{N\setminus M}^{h}$ is a subfunctor of $G_{N\setminus M}^{r}$.  The
corresponding morphisms of schemes, $G^{h}_{N}\rightarrow G^{r}_{N}$
and $G_{N\setminus M}^{h}\rightarrow G_{N\setminus M}^{r}$, are closed
embeddings.  To see this, observe that $G^{h}_{N}$ is defined locally
by the vanishing of the coordinates $\gamma ^{x}_{b}$ on $G^{r}_{N}$
with $x\in N_{a}$, $b\in N_{c}$, for $a\not =c$. \qed 
\end{remark}

We will now prove the two theorems stated at the beginning of this section.

\medskip

\noindent {\em Proof of Theorem \ref{t:Hilb-scheme-I}:} We shall apply
Proposition \ref{p:relative-representability} to represent $H_{T}^{h}$
in $G^{h}_{N\setminus M}$.

\nobreak

{\em Step 1: $H_{T}^{h}$ is a Zariski sheaf.}  Let
$f_{1},\ldots,f_{k}$ generate the unit ideal in $R$.  To give a
homogeneous submodule $L\subseteq {R\otimes T}$, it is equivalent to
give a compatible system of homogeneous submodules $L_{i}\subseteq
{R_{f_{i}}\otimes T}$.  The homogeneous component $L_{a}$ is locally
free of rank $h(a)$ if and only if the same holds for each
$(L_{i})_{a}$.

{\em Step 2: For all $R\in \kAlgCateg $ and $L\in H_{T}^{h}(R)$, $M$
generates $({R\otimes T})/L$ as an $R$-module.}  Localizing at each
$P\in \Spec R$, it suffices to prove this when $(R,P)$ is a local
ring.  Then for all $a\in A$, the $R$-module $({R \otimes T_{a}})
/L_{a}$ is free of finite rank $h(a)$.  By Nakayama's Lemma, $RM_{a} =
({R\otimes T_{a}})/L_{a}$ if and only if $K M_{a} = ({K \otimes
T_{a}})/L_{a}$, where $K = R/P$ is the residue field.  The latter
holds by hypothesis (iii).

{\em Step 3: We have a canonical natural transformation $\eta \colon
H_{T}^{h}\rightarrow G^{h}_{N\setminus M}$}.  It follows from Step 2
that the canonical homomorphism ${R\otimes N}\rightarrow ({R\otimes
T})/L$ is surjective.  If $L'$ denotes its kernel, it further follows
that $M$ generates $({R\otimes N})/L'$.  Hence we have $L'\in
G^{h}_{N\setminus M}(R)$, and the rule $\eta _{R}(L) = L'$ clearly
defines a natural transformation.  Note that $G^{h}_{N\setminus M}$
makes sense as a scheme functor by hypothesis (i).

{\em Step 4: The functors $\eta ^{-1}G^{h}_{N\setminus B}$ are
represented by affine schemes.}  Let $B\subseteq M$ be any homogeneous
subset with $|B_{a}| = h(a)$ for all $a\in A$, so $G^{h}_{N\setminus
B}$ is a standard affine chart in $G^{h}_{N\setminus M}$.  In
functorial terms, $G^{h}_{N\setminus B}(R)$ describes quotients
$({R\otimes N})/L'$ that are free with basis $B$.  Hence $\eta
^{-1}G^{h}_{N\setminus B}(R)$ consists of those $L\in H_{T}^{h}(R)$
such that $({R\otimes T})/L$ is free with basis $B$. Given such an
$L$, we define coordinates $\gamma^{x}_{b}\in R$ for all $a\in A$ and
all $x\in T_{a}$, $b\in B_{a}$ by requiring that $x - \sum_{b \in B_a}
{\gamma^{x}_{b} \cdot b }$ is in $ L$.  For $x\in N$, the coordinates
$\gamma ^{x}_{b}$ of $L$ coincide with the Grassmann functor
coordinates of $\eta _{R}(L)$, so there is no ambiguity of notation.  In
particular, they satisfy
\begin{equation}\label{e:delta-condition}
\gamma ^{x}_{b} =  \delta _{x,b}\quad \text{for $x\in B$}.
\end{equation}
They also clearly satisfy a syzygy condition similar to
\eqref{e:J-in-L-condition}, for every linear relation $\sum
_{x}{c_{x}\cdot x} = 0$, $c_{x}\in k$, holding among elements $x\in
T_{a}$.  Namely,
\begin{equation}\label{e:syzygy-condition}
\sum _{x}c_{x}\cdot \gamma ^{x}_{b} = 0\quad \text{for all $a\in A$,
$b\in B_{a}$.}
\end{equation}
Finally, since $L$ is an $F$-submodule, the coordinates $\gamma
^{x}_{b}$ satisfy
\begin{equation}\label{e:F-condition}
\gamma ^{fx}_{b} = \sum _{b'\in B_{a}}\gamma ^{x}_{b'}\gamma
^{fb'}_{b}\quad \text{for all $a,c\in A$ and all $x\in T_{a}$, $f\in
F_{ac}$, $b\in B_{c}$}.
\end{equation}

Conversely, suppose we are given elements $\gamma ^{x}_{b}\in R$
satisfying equations \eqref{e:delta-condition}--\eqref{e:F-condition}.
We fix attention on an individual degree $a$ for the moment.
The elements $\gamma ^{x}_{b}$ for $x\in T_{a}$, $b\in B_{a}$ can be
viewed as the entries of a (typically infinite) matrix defining a
homomorphism of free $R$-modules
\begin{equation}\label{e:Ta->Ba}
\phi_{a} \colon R^{T_{a}}\rightarrow R^{B_{a}}.
\end{equation}
Equation \eqref{e:syzygy-condition} ensures that $\phi_{a} $ factors
through the canonical map $R^{T_{a}}\rightarrow {R\otimes T_{a}}$,
inducing $\phi'_{a}\colon {R\otimes T_{a}}\rightarrow R^{B_{a}}$.
Equation \eqref{e:delta-condition} ensures that $\phi'_{a}$ is the
identity on $B_{a}$.  Let $L_{a}$ be the kernel of $\phi'_{a}$. We
conclude that $({R\otimes T_{a}})/L_{a}$ is free with basis $B_{a}$.
Considering all degrees again, equation \eqref{e:F-condition} ensures
that the homogeneous $R$-submodule $L\subseteq {R\otimes T}$ thus
defined is an $F$-submodule.  We have given correspondences in both
directions between elements $L\in \eta ^{-1}G^{h}_{N\setminus B}(R)$
and systems of elements $\gamma ^{x}_{b}\in R$ satisfying
\eqref{e:delta-condition}--\eqref{e:F-condition}.  These two
correspondences are mutually inverse and natural in $R$.  By \cite[\S
I.4]{EH}, this shows that $\eta ^{-1}G^{h}_{N\setminus B}$ is
represented by an affine scheme over $k$.

{\em Step 5.}  It now follows from Proposition
\ref{p:relative-representability} that $H_{T}^{h}$ is represented by a
scheme over $G^{h}_{N\setminus M}$, the morphism $H_{T}^{h}\rightarrow
G^{h}_{N\setminus M}$ being given by the natural transformation $\eta$
from Step 3.  Up to this point, we have only used hypotheses (i) and
(iii).

{\em Step 6.  The morphism corresponding to $\eta \colon
H_{T}^{h}\rightarrow G^{h}_{N\setminus M}$ is a closed embedding.}  It
is enough to prove this locally for the restriction of $\eta $ to the
preimage of $G^{h}_{N\setminus B}$.  This restriction corresponds to
the morphism of affine schemes given by identifying the coordinates
$\gamma ^{x}_{b}$ on $G^{h}_{N\setminus B}$ with those of the same
name on $\eta ^{-1}G^{h}_{N\setminus B}$.  To show that it is a closed
embedding, we must show that the corresponding ring homomorphism is
surjective.  In other words, we claim that the elements $\gamma
^{x}_{b}$ with $x\in N$ generate the algebra $k[\{\gamma
^{x}_{b}\}]/I$, where $I$ is the ideal generated by
\eqref{e:delta-condition}--\eqref{e:F-condition}.  Consider the
subalgebra generated by the $\gamma ^{x}_{b}$ with $x\in N$. Let $g
\in G$.  If $\gamma ^{x}_{b}$ belongs to the subalgebra for all $b\in
B$, then so does $\gamma ^{gx}_{b}$, by equation \eqref{e:F-condition}
and hypothesis (iv).  Since $G$ generates $F$, and $N$ generates $T$
as an $F$-module by hypothesis (ii), we conclude that $\gamma
^{x}_{b}$ lies in the subalgebra for all $x$.  Theorem
\ref{t:Hilb-scheme-I} is now proved. \qed

\medskip

A description of the Hilbert scheme in terms of affine charts is
implicit in the proof above.  There is a chart for each homogeneous
subset $B$ of $M$ with $h(a)$ elements in each degree $a$, and the
coordinates on that chart are the $\gamma ^{x}_{b}$ for homogeneous
elements $x$ generating $N$.  
Local equations are derived
from \eqref{e:delta-condition}--\eqref{e:F-condition}.

\medskip

\noindent {\em Proof of Theorem \ref{t:Hilb-scheme-II}:} We will show
that Proposition \ref{p:closed-condition} applies to $H_{T}^{h}
\rightarrow H_{T_{D}}^{h}$.

\nobreak 

{\em Step 1: For $L_{D}\in H_{T_{D}}^{h}(R)$, let $L'\subseteq
{R\otimes T}$ be the $F$-submodule generated by $L_{D}$.  Then the
$R$-module $({R\otimes T_{a}})/L'_{a}$ is finitely generated in each
degree $a\in A$. } Take $E$ as in (v) and let $Y$ be a finite
generating set of the $k$-module $T_{a}/\sum _{b\in D} E_{ba}(T_{b})$.
Since $E$ is finite, the sum can be taken over $b$ in a finite set of
degrees $D'\subseteq D$.

For $b\in D'$, the $R$-module $({R\otimes T_{b}})/L'_{b}$ is locally
free of rank $h(b)$, and hence generated by a finite set $M_{b}$.  For
all $x\in {R\otimes T_{b}}$ there exist coefficients $\gamma
^{x}_{v}\in R$ (not necessarily unique, as $M_{b}$ need not be a
basis) such that $x \equiv \sum _{v\in M_{b}}{\gamma ^{x}_{v} \cdot v}
\pmod {L'_{b}}$.  For all $g\in E_{ba}$ we have $gx \equiv \sum _{v\in
M_{b}}{\gamma ^{x}_{v} \cdot gv} \pmod{L'_{a}}$.  This shows that the
finite set $Z =\bigcup _{b\in D', g\in E_{ba}} g(M_{b})$ generates the
image of $R\otimes \sum _{b\in D} E_{ba} (T_{b})$ in $({R \otimes
T_{a}})/L_{a}'$, and therefore $Y\cup Z$ generates $({R \otimes
T_{a}})/L'_{a}$.

{\em Step 2: $H_{T}^{h}$ is a subfunctor of $H_{T_{D}}^{h}$.}
Equivalently, for all $k$-algebras $R$, the map
$H_{T}^{h}(R)\rightarrow H_{T_{D}}^{h}(R)$, $L \mapsto L_D$ is
injective.  We will prove that if $L'\subseteq {R\otimes T}$ is the
$F$-submodule generated by $L_{D}$, then $L' = L$.  Localizing at a
point $P\in \Spec R$, we can assume that $(R,P)$ is local, and hence
the locally free modules $({R\otimes T_{a}})/L_{a}$ are free.  Fix a
degree $a\in A$, and let $B_{a}$ be a free module basis of $({R\otimes
T_{a}})/L_{a}$.  Then $B_{a}$ is also a vector space basis of
$({K\otimes T_{a}})/({K \otimes L_{a}})$, where $K = R/P$ is the
residue field.  In particular, $\dim ({K\otimes T_{a}})/({K \otimes
L_{a}}) = |B_{a}| = h(a)$.  By (vi) we have $\dim ({K\otimes
T_{a}})/({K \cdot L'_{a}}) \leq h(a)$, and hence ${K \cdot L'_{a}} =
{K\otimes L_{a}}$, since $L' \subseteq L$.  By Step 1, the $R$-module
$({R\otimes T_{a}})/L'_{a}$ is finitely generated, so Nakayama's Lemma
implies that $B_{a}$ generates $({R \otimes T_{a}})/L'_{a}$.  Since
$B_{a}$ is independent modulo $L_{a} \supseteq L'_{a}$, it follows
that $L'_{a} = L_{a}$.

{\em Step 3: The condition that $({S\otimes T_{a}})/L'_{a}$ be locally
free of rank $h(a)$ is closed.}  More precisely, fix a $k$-algebra $R$
and $L_D \in H_{T_{D}}^{h}(R)$.  Given an $R$-algebra $\phi \colon
R\rightarrow S$, let $L'\subseteq {S\otimes T}$ be the $F$-submodule
generated by $H_{T_{D}}^{h}(\phi )L_{D} = {S\otimes _{R}L_{D}}$.  Then
the condition that $({S\otimes T_{a}})/L'_{a}$ be locally free of rank
$h(a)$ is a closed condition on $S$.
To see this, let $L^{0}$ be the $F$-submodule of ${R\otimes T}$
generated by $L_{D}$, that is, the $L'$ for the case $S = R$.  By Step
1, $({R\otimes T_{a}})/L^{0}_{a}$ is finitely generated, say by a set
$X$.  By (vi) and Nakayama's Lemma, $({R_{P}\otimes
T_{a}})/(L^{0}_{a})_{P}$ is generated by at most $h(a)$ elements of
$X$, for every $P\in \Spec R$.  For every subset $B\subseteq X$ with
$|B|=h(a)$ elements, the set of points $P\in \Spec R$ where $B$
generates $({R_{P}\otimes T_{a}})/(L^{0}_{a})_{P}$ is an open set
$U_{B}$, and these open sets cover $\Spec R$.  The property that a
condition on $R$-algebras is closed is local with respect to the base
$R$.  Therefore, replacing $R$ by some localization $R_{f}$, we can
assume that a single set $B$ with $h(a)$ elements generates
$({R\otimes T_{a}})/L^{0}_{a}$.  Then $B$ also generates $({S\otimes
T_{a}})/L'_{a}$ for every $R$-algebra $S$.

A presentation of the $S$-module $({S\otimes T_{a}})/L'_{a} = S\otimes
_{R}({(R\otimes T_{a})/L^{0}_{a}})$ is given by the generating set
$B$, modulo those relations on $b \in B$ that hold in $({R\otimes
T_{a}})/L^{0}_{a}$:
\begin{equation}\label{e:syzygies-mod-L0}
\sum _{b\in B} c_{b} \cdot b \equiv 0 \pmod{L^{0}_{a}}, \quad c_{b}\in
R.
\end{equation}
Thus $({S\otimes T_{a}})/L'_{a}$ is locally free of rank $h(a)$ if and
only if it is free with basis $B$, if and only if all coefficients
$c_{b}$ of all syzygies in \eqref{e:syzygies-mod-L0} vanish in $S$,
{\it i.e.}, $\phi (c_{b}) = 0$.  This condition is closed, with
defining ideal $I\subseteq R$ generated by all the coefficients
$c_{b}$.

{\em Step 4: The subfunctor $H_{T}^{h}\rightarrow H_{T_{D}}^{h}$ is
represented by a closed subscheme.}  By Step 2, $H_{T}^{h}$ is a
subfunctor, and by Step 1 in the proof of
Theorem~\ref{t:Hilb-scheme-I}, it is a Zariski sheaf.  In Step 2 we
saw that $L_{D}\in H_{T_{D}}^{h}(S)$ is in the image of $H_{T}^{h}(S)$
if and only if the $F$-submodule $L'$ it generates belongs to
$H_{T}^{h}(S)$.  By Step 3, this is a closed condition, since it is
the conjunction of the conditions that $({S \otimes T_{a}})/L'_{a}$ be
locally free of rank $h(a)$, for all  $a\in A$.
Theorem~\ref{t:Hilb-scheme-II} now follows from Proposition
\ref{p:closed-condition}.  \qed

\smallskip

The algorithmic problem arising from Theorem~\ref{t:Hilb-scheme-II} is
to give equations on $H_{T_{D}}^{h}$ which define the closed subscheme
$H_{T}^{h}$.  We assume that we already have a description of an
affine open subset $U\subseteq H_{T_{D}}^{h}$ as $\Spec R$ for some
$k$-algebra $R$ (see the paragraph following the proof of
Theorem~\ref{t:Hilb-scheme-I} above). The embedding of $U=\Spec R$
into $H_{T_{D}}^{h}$ corresponds to a {\it universal element} $L\in
H_{T_{D}}^{h}(R)$.  The ideal $I\subseteq R$ defining the closed
subscheme $H_{T}^{h}\cap U$ is generated by separate contributions
from each degree $a$, determined as follows.  Construct the finite set
$X=Y\cap Z\subseteq T_{a}$ in Step 1, and compute the syzygies of $X$
modulo $L^{0}_{a}$, where $L^{0}\subseteq {R\otimes T}$ is the
$F$-submodule generated by $L$.  These syzygies are represented by the
(perhaps infinitely many) rows of a matrix $\Gamma $ over $R$, with
columns indexed by the finite set $X$.  The content of hypothesis (vi)
is that the minors of size $|X|-h(a)$ in $\Gamma $ generate the unit
ideal in $R$.  The contribution to $I$ from degree $a$ is the {\em
Fitting ideal} $\, I_{|X|-h(a)+1}(\Gamma ) \,$ generated by the minors
of size $|X|-h(a)+1$.  In fact, the vanishing of these minors,
together with the fact that $I_{|X|-h(a)}(\Gamma )$ is the unit ideal,
is precisely the condition that the submodule $L'_{a}\subseteq
{R\otimes T_{a}}$ generated by the rows of $\Gamma $ should have
$({R\otimes T_{a}})/L'_{a}$ locally free of rank $h(a)$.

If $k$ is Noetherian, so $H_{T_{D}}^{h}$ is a Noetherian scheme, then
$H_{T}^{h}$ must be cut out as a closed subscheme by the equations
coming from a finite subset $E\subseteq A$ of the degrees.  As we
shall see, this is also true when $T$ is a multigraded polynomial
ring, even if the base ring $k$ is not Noetherian.  Finding such a set
$E$ amounts to finding an isomorphism $H_{T}^{h}\cong H_{T_{E}}^{h}$.
Satisfactory choices of $D$ and $E$ for multigraded Hilbert schemes
will be discussed in the next section.

\smallskip

Here is a simple example, taken from \cite[\S VI.1]{Bay}, to
illustrate our results so far.

\begin{example} \label{daveRunning}
Let $A = \{3,4\}$, $T_3 \simeq k^4$ with basis $\{x^3,x^2 y, x y^2,
y^3\}$, $T_4 \simeq k^5$ with basis $\{x^4,x^3 y, x^2 y^2, x y^3,
y^4\}$, and $F = F_{3,4} = \{x,y\}$, {\it i.e.}, the operators are
multiplication by variables.  Fix  $h(3) = h(4) =
1$, and $D = \{3\}$. Then $H_{T_{D}}^{h}$ is the projective space
$\P^3$ parametrizing rank $1$ quotients of $T_3$, where 
$({c_{123} : c_{124} : c_{134} : c_{234}})  \in \fun{\P^3}(R) $
corresponds to the $R$-module $L_D = L_{3}$ generated by the $2 \times
2$-minors of 
\begin{equation}\label{e:bayer-coords}
\left(
\begin{array}{cccc}
c_{234} & -c_{134} & c_{124} & - c_{123} \\
x^3 & x^2 y & x y^2 & y^3 
\end{array} \right).
\end{equation}
The Hilbert scheme $H_{T}^{h}$ is the projective line $\P^1$ embedded
as the twisted cubic curve in $H_{T_{D}}^{h} \simeq \P^3$ defined by
the quadratic equations
\begin{equation}\label{e:bayer-example}
c_{134} c_{124} - c_{123} c_{234} = c_{124}^2 - c_{123} c_{134} =
c_{134}^2 - c_{124} c_{234} = 0 .
\end{equation}
\end{example}

\section{Constructing the multigraded Hilbert scheme}
\label{sec3}

We now take up our primary application of
Theorems~\ref{t:Hilb-scheme-I} and \ref{t:Hilb-scheme-II}, the
construction of multigraded Hilbert schemes.  Let $S = k[\x ] =
k[x_{1},\ldots,x_{n}]$ be a polynomial ring over $k$, with a
multigrading $S = \bigoplus _{a}S_{a}$ induced by a degree function
$\deg \colon \N ^{n}\rightarrow A$, with $\deg (x^{u}) = \deg (u)$, as
in the introduction.  Here $A$ is an abelian group, or the
subsemigroup $A_+$ generated by $\deg(x_i) = a_i$ for $1=1,\ldots,n$.
As our $k$-module with operators $(T,F)$ we take $T=S$, with $F$ the
set of all multiplications by monomials.  More precisely, $F_{ab}$
consists of multiplications by monomials of degree $b-a$, for all
$a,b\in A$.  Then an $F$-submodule $L\subseteq {R \otimes S}$ is an
ideal of ${R \otimes S} = R[\x ]$ which is homogeneous with
respect to the $A$-grading.

Fix a Hilbert function $h\colon A\rightarrow \N $, and let $H_{S}^{h}$
be the Hilbert functor. For any $k$-algebra $R$, the set
$H_{S}^{h}(R)$ consists of admissible homogeneous ideals $L\subseteq
{R \otimes S}$ with Hilbert function $h$.  Theorem \ref{t:multi-Hilb}
states that the functor $H_{S}^{h}$ is represented by a
quasiprojective scheme.  For the proof we need two 
combinatorial results.

\begin{lemma}[Maclagan \cite{M}]\label{l:Diane}
Let $C$ be a set of monomial ideals in $k[\x ]$ which is an antichain,
that is, no ideal in $C$ contains another.  Then $C$ is
finite. \end{lemma}

Let $I \subseteq k[\x ]$ be a monomial ideal and $\deg \colon \N
^{n}\rightarrow A$ a multigrading.  The monomials not in $I$ are
called the {\em standard monomials} for $I$.  The value $h_I(a)$ of
the Hilbert function $h_I$ at $a \in A$ is the number of standard
monomials in degree $a$.

\begin{proposition}\label{p:Diane-prime}
Given a multigrading $\deg \colon \N^{n}\rightarrow A$ and a Hilbert
function $h\colon A\rightarrow \N $, there is a finite set of degrees
$D\subseteq A$ with the following two properties:
\begin{itemize}
\item [(g)] Every monomial ideal with Hilbert function $h$ is
generated by monomials of degree belonging to $D$, and
\item [(h)] every monomial ideal $I$ generated in degrees $D$
satisfies: if $h_I(a) = h(a)$ for all $a \in D$, then $h_I(a) = h(a)$
for all $a \in A$.
\end{itemize}
\end{proposition}

Our labels for these properties are mnemonics for {\em
generators} and {\it Hilbert function}.

\begin{proof}
Let $C$ be the set of all monomial ideals with Hilbert function $h$.
By Lemma \ref{l:Diane}, $C$ is finite.  Let $D_{0}$ be the set of all
degrees of generators of ideals in $C$.  Now let $C_{0}$ be the set of
monomial ideals that are generated in degrees in $D_{0}$ and whose
Hilbert function agrees with $h$ on $D_{0}$.  By Lemma \ref{l:Diane}
again, $C_{0}$ is finite.  If $C_{0} = C$, then $D_{0}$ is the
required $D$.  Otherwise, for each ideal $I\in C_{0} \backslash C$, we
can find a degree $a$ with $h_I(a) \not= h(a)$. Adjoining finitely
many such degrees to $D_{0}$, we obtain a set $D_{1}$ such that every
monomial ideal generated in degrees $D_{0}$ and having Hilbert
function $h$ in degrees $D_{1}$ belongs to $C$.  Now we define $C_{1}$
in terms of $D_{1}$ as we defined $C_{0}$ in terms of $D_{0}$, namely,
$C_{1}$ is the set of monomial ideals generated in degrees $D_{1}$ and
with Hilbert function $h$ on $D_{1}$.  By construction, we have
$C_{1}\cap C_{0} = C$.  Iterating this process, we get a sequence
$C_{0}, C_{1}, C_{2},\ldots $ of sets of monomial ideals with $C_{i}
\cap C_{i+1} = C$ for all $i$, and finite sets of degrees $D_{0}
\subseteq D_{1} \subseteq D_{2} \subseteq \cdots$.  Here $D_i$ are the
degrees such that every monomial ideal generated in degrees $D_{i-1}$
and with Hilbert function $h$ in degrees $D_{i}$ belongs to $C$, and
$C_i$ are the monomial ideals generated in degrees $D_{i}$ and with
Hilbert function $h$ on $D_{i}$.  We claim that this sequence
terminates with $C_{k} = C$ for some $k$.

Given an ideal $I_j \in C_j$, its {\it ancestor} in $C_{i}$ for $i<j$
is the ideal $I_{i}$ generated by the elements of degrees $D_{i}$ in
$I_{j}$.  We say that $I_j$ is a {\em descendant} of its ancestors.
If $I_{j}$ is a descendant of $I_{i}$, then $I_i \subseteq I_j$, and
$I_i \in C $ implies $I_{i}=I_j$.  Suppose, contrary to our claim,
that $C_{k}\not =C$ for all $k$.  Since $C_{0}$ is finite, there is an
$I_{0}\in C_{0} \backslash C$ with descendants in $C_{k} \backslash C
$ for infinitely many $k$, and hence for all $k>0$.  Among its
descendants in $C_{1}$ must be one, call it $I_{1}$, with descendants
in $C_{k} \backslash C$ for all $k>1$.  Iterating, we construct a
sequence $I_{0},I_{1},\ldots$ with $I_{k}\in C_{k}$ and $I_{k+1}$ a
descendant of $I_{k}$. By the ascending chain condition, $I_{k} =
I_{k+1}$ for some $k$.  But this implies $I_{k}\in C$, a
contradiction.  We conclude that $C_k = C$ for some $k$, and $D =
D_{k}$ is the required set of degrees.
\end{proof}

\begin{lemma}\label{l:E}
Given a multigrading $\deg \colon \N ^{n}\rightarrow A$, let
$D\subseteq A$ be a subset of the degrees and $J = \langle x^{u}:\deg (u)\in
D \rangle$ the ideal generated by all monomials with degree in $D$.  If $a\in
A$ is a degree such that $h_{J}(a)$ is finite, then there is a finite
set of monomials $E\subseteq \bigcup _{b\in D} F_{ba}$ such that
$S_{a}/\sum _{b\in D} E_{ba}(S_{b})$ is finitely generated.
\end{lemma}

\begin{proof}
Choose an expression for each minimal monomial in $J_{a}$ as
$x^{v}x^{u}$ for some $x^{u}\in S_{b}$, $b\in D$, and let $E$ be the
set of monomials $x^{v}$ that occur.  For all $x^{r}\in J_{a}$, we
have $x^{r} = x^{q}x^{v}x^{u}$ for some minimal $x^{v}x^{u}\in J_{a}$,
and $\deg (q) = 0$.  Hence $x^{r} = x^{v}(x^{q}x^{u})\in
E_{ba}(S_{b})$.  This shows that set of all standard monomials of
degree $a$ for $J$ spans $S_{a}/\sum _{b\in D} E_{ba}(S_{b})$.  This
set is finite, by hypothesis.
\end{proof}

We are now ready to construct the multigraded Hilbert scheme.
In our proof, the condition (h) in Proposition \ref{p:Diane-prime} 
will be replaced by the following weaker condition.
\begin{itemize}
\item [(h$'$)] every monomial ideal $I$ generated in degrees $D$
satisfies: if $h_I(a) = h(a)$ for all $a \in D$, then $h_I(a) \leq
h(a)$ for all $a \in A$.
\end{itemize}
Proposition \ref{p:Diane-prime} holds verbatim for ``(g) and (h$'$)''
instead of ``(g) and (h)''.  We fix a {\em term order} on $\N^n$, so
that each ideal $L\subseteq K[\x ]$, with $K \in \kAlgCateg $ a field,
has an {\em initial monomial ideal} $\init (L)$.  The Hilbert function
of $\init(L) $ equals that of $L$.

\medskip

\noindent {\em Proof of Theorem~\ref{t:multi-Hilb}:} By definition,
$F$ is the system of operators on $S = k[\x]$ given by multiplication
by monomials.  We first verify the hypotheses of
Theorem~\ref{t:Hilb-scheme-I} for $(S_{D},F_{D})$, where $D\subseteq
A$ is any finite subset of the degrees.  Let $C$ be the set of
monomial ideals generated by elements of degrees in $D$, and with
Hilbert function agreeing with $h$ on $D$. By Lemma \ref{l:Diane}, the
set $C$ is finite.  Let $M'$ be the union over all $I\in C$ of the set
of standard monomials for $I$ in degrees $D$.  Then $M'$ is a finite
set of monomials which spans the free $k$-module $S_{D}/I_{D}$ for all
$I\in C$.

The monomials of degree zero in $S$ form a finitely generated
semigroup.  Let $G'_{0}$ be a finite generating set for this
semigroup, so that $S_0$ is the $k$-algebra generated by $G'_{0}$.
Every component $S_{a}$ is a finitely generated $S_{0}$-module.
For each $a\in A$, let $F'_{a}$ be a finite set of monomials
generating $S_{a}$ as an $S_{0}$-module.  Then every monomial of degree
$a$ is the product of a monomial in $F'_{a}$ and zero or more
monomials in $G'_{0}$.  For $b, c\in D$, let $G_{bc} \subseteq F_{bc}$
consist of multiplications by monomials in $F'_{c-b}$, if $b\not =c$,
or in $G'_{0}$, if $b=c$.  Then $G = \bigcup _{b,c} G_{bc}$ is finite
and generates $F_{D}$ as a category.

Our construction is based on the following finite set of monomials:
\begin{equation}\label{ChoiceOfN}
N' = GM' \cup \bigcup _{a\in D} F'_{a} .
\end{equation}
Setting $M = kM'$, $N = kN'$, it is obvious that $M$, $N$ and $G$
satisfy hypotheses (i), (ii) and (iv) of Theorem
\ref{t:Hilb-scheme-I}.  For (iii), fix a field $K\in \kAlgCateg $ and
an element $L_{D}\in H_{S_{D}}^{h}(K)$.  Let $L \subseteq {K\otimes
S}$ be the ideal generated by $L_{D}$ and $I$ the monomial ideal
generated by $\init (L)_{D}$.  Equivalently, $I$ is the ideal
generated by the leading monomials of elements of $L_D$. Therefore $I$
belongs to $C$ and its standard monomials span $({K\otimes
S_{D}})/L_{D}$.  We conclude that $M'$ spans $({K\otimes
S_{D}})/L_{D}$, which proves (iii).  We have now shown that
$H^{h}_{S_{D}}$ is represented by a quasiprojective scheme for every
finite set of degrees $D$.

It remains to verify hypotheses (v) and (vi) of Theorem
\ref{t:Hilb-scheme-II} for a suitable choice of $D$.  Let $D$ be any
finite subset of $A$ that satisfies the conditions (g) and (h$'$).  We
assume that there exists a monomial ideal $I$ generated in degrees $D$
and satisfying $h_{I}(a) = h(a)$ for all $a\in D$.  Otherwise the
Hilbert functor and Hilbert scheme are empty, so the result holds
vacuously.  By condition (h$'$), $h_{I}(a)$ is finite for all $a\in
A$.  The ideal $J$ in Lemma~\ref{l:E} contains $I$, so $h_{J}(a)$ is
also finite.  For hypothesis (v), we can therefore take $E$ as given
by Lemma~\ref{l:E}.

For (vi), we fix $K$ and $L_{D}\in H_{S_{D}}^{h}(K)$ as we did for
(iii), and again let $L$ be the ideal generated by $L_{D}$ and $I$ the
ideal generated by $\init (L)_{D}$.  Our assumption on $D$ implies
that the Hilbert function of $I$ satisfies $h_I(a) \leq h(a)$ for all
$a \in A$.  Since $I\subseteq \init (L)$ it follows that $h_{L}(a) =
h_{\init (L)}(a) \leq h_I(a) \leq h(a)$ for all $a \in A$.  This
establishes hypothesis (vi). We have proved that the Hilbert functor
$H_{S}^{h}$ is represented by a closed subscheme of $H_{S_{D}}^{h}$.
\qed

\medskip

Our ultimate goal is to compute the scheme $H_{S}^{h}$ effectively.
One key issue is to identify suitable finite sets of degrees.  A
subset $D$ of the abelian group $A$ is called {\em supportive} for $h$
if $D$ is finite and the conditions (g) and (h$'$) are satisfied.  The
last two paragraphs in the proof of Theorem~\ref{t:multi-Hilb} establish
the following result.

\begin{corollary}\label{p:multi-closed}
Take $S$ and $h \colon A \rightarrow \N$ as in
Theorem~\ref{t:multi-Hilb}.  If the set of degrees $D \subseteq A$ is
supportive then the canonical morphism $H_{S}^{h}\rightarrow
H_{S_{D}}^{h}$ is a closed embedding.
\end{corollary}

\begin{remark}\label{r:props} Corollary~\ref{p:projective} follows
immediately from this result and Remark~\ref{r:projective}.
Using Remark~\ref{r:refinement}, Proposition~\ref{p:refinement} also
follows.
\end{remark}

Consider one further condition on Hilbert functions and subsets of
degrees:
\begin{itemize}
\item [(s)] For every monomial ideal $I$ with $h_{I} = h$, the syzygy
module of $I$ is generated by syzygies $x^{u}x^{v'} = x^{v}x^{u'} =
\lcm (x^{u},x^{v})$ among generators $x^{u}$, $x^{v}$ of $I$ such that
$\deg \lcm (x^{u},x^{v})\in D$ ({\it i.e.}, all minimal S-pairs have
their degree in $D$).
\end{itemize}
A subset $D$ of $A$ is called {\em very supportive} for a given
Hilbert function $h \colon A \rightarrow \N$ if $D$ is finite and the
conditions (g), (h) and (s) are satisfied.  It follows from
Proposition \ref{p:Diane-prime} that a very supportive set of degrees
always exists.

\begin{theorem} \label{t:multi-iso}
Take $S$ and $h \colon A \rightarrow \N$ as in
Theorem~\ref{t:multi-Hilb}.  If the set of degrees $D \subseteq A$ is
very supportive then the canonical morphism $H_{S}^{h}\rightarrow
H_{S_{D}}^{h}$ is an isomorphism.
\end{theorem}

\begin{example} 
Let $S = k[x,y,z]$ with the $\Z$-grading $\deg(x) = \deg(y) = 1$ and $
\deg(z) = -1$ and fix the Hilbert function $ h(a) = 2 $ for all $a \in
\Z$.  This example is typical in that both the support of $h$ and the
set of monomials in any fixed degree are infinite.  There are eight
monomial ideals with this Hilbert function:
\begin{eqnarray*}
  \langle x^2 z^2, y \rangle , 
  \langle x^2, y z \rangle , 
  \langle x^2 z, x y, y z \rangle , 
  \langle x^2 z, y^2, y z \rangle , \\
  \langle y^2 z^2, x \rangle , 
  \langle y^2, x z \rangle , 
  \langle y^2 z, x y, x z \rangle , 
  \langle y^2 z, x^2, x z \rangle.
\end{eqnarray*}
The set $D = \{0,1,2\}$ is very supportive, so the Hilbert
scheme $H^h_{S}$ is isomorphic to $H^h_{S_D}$.  It 
can be checked that this scheme is smooth of
dimension $4$ over $\Spec k$.   \qed
\end{example}

For the proof of Theorem \ref{t:multi-iso}
we need a variant of Gr\"obner bases for ideals in the
polynomial ring over a local ring $R$.  Let $(R,P)$ be a local ring
satisfying
\begin{equation}\label{e:orderly}
\bigcap _{m} P^{m} =0.
\end{equation}
This holds for example if $R$ is complete or Noetherian.  Let $R[\x ]
= R[x_{1},\ldots,x_{n}]$ and fix a term order on $\N ^{n}$.  This
induces a lexicographic order $<$ on the set $(-\N ) \times \N ^{n}$,
in which $(-d,e)<(-d',e')$ if $-d<-d'$ or if $d=d'$ and $e<e'$ in the
given term order.  The lexicographic order is not well-ordered, but
has the property that if
\[
(-d_{1},e_{1})>(-d_{2},e_{2})>\cdots 
\]
is an infinite strictly decreasing chain, then the sequence
$d_{1},d_{2},\ldots$ is unbounded.

\begin{definition}\label{d:ord-in}
The {\em order} $\ord (a)$ of a nonzero element $a\in R$ is the unique
integer $m$ such that $a\in P^{m}\backslash P^{m+1}$, which exists by
\eqref{e:orderly}.  The {\em initial term} $\init (p)$ of a nonzero
polynomial $p\in R[\x ]$ is the term $ax^{e}$ of $p$ for which the
pair $(-\ord (a), e)\in (-\N ) \times \N ^{n}$ is lexicographically
greatest.
\end{definition}

The definition of initial term is compatible with the following
filtration of $R[\x ]$ by $R$-submodules: given $(-d,e)\in (-\N )
\times \N ^{n}$, we define $R[\x ]_{\leq (-d,e)}$ to be the set of
polynomials $p$ such that for every term $bx^{h}$ of $p$, we have
$(-\ord(b),h)\leq (-d,e)$.  We also define $R[\x ]_{<(-d,e)}$ in the
obvious analogous way.  Then $\init (p) = ax^{e}$, with $\ord (a) =
d$, if and only if $p\in R[\x ]_{\leq (-d,e)}\backslash R[\x
]_{<(-d,e)}$.

Consider a set of nonzero polynomials $F\subseteq R[\x ]$ satisfying
the restriction:
\begin{equation}\label{e:restriction}
\text{for all $f\in F$, the initial term of $f$ has coefficient $1$.}
\end{equation}
Let $I = \langle F \rangle$ be the ideal in $R[\x ]$ generated by $F$.

\begin{definition}\label{d:Grobner}
A set $F$ satisfying \eqref{e:restriction} is a {\em Gr\"obner basis}
of the ideal $ I = \langle F \rangle $ if for all nonzero $p\in I$,
the initial term $\init (p)$ belongs to the monomial ideal generated
by the set of initial terms $\init (F) = \{\init (f):f\in F \}$.
\end{definition}

In general, we do not have $\init (pq) = \init (p)\init (q)$, but this
does hold when $\init (p) = ax^{e}$, $\init (q) = bx^{h}$ with $\ord
(ab) = \ord (a) + \ord (b)$.  Thus condition \eqref{e:restriction}
implies $\init (pf) = \init (p)\init (f)$, and in particular
\begin{equation}\label{e:in(axf)}
\init (ax^{e}f) = ax^{e}\init (f) \qquad
\hbox{for $f\in F$}.
\end{equation}
Without this condition, even a one-element set $F$
could fail to be a Gr\"obner basis.

If $F$ is a Gr\"obner basis, the standard monomials for the monomial
ideal $ \langle \init (F) \rangle$ are $R$-linearly independent modulo
$I$, since every nonzero element of $I$ has an initial term belonging
to $\langle \init (F) \rangle$.  There is a reformulation of the
Gr\"obner basis property in terms of a suitably defined notion of
$F$-reducibility.

\begin{definition}\label{d:reducible}
A polynomial $p\in R[\x ]$ is {\em $F$-reducible} if $p=0$ or if
$\init (p)= ax^{e}$, $d = \ord (a)$ and for all $m\geq 0$ there exists
an expression
\begin{equation}\label{e:reducible}
p \equiv  \sum b_{i}x^{h_{i}}f_{i} \quad \pmod {P^{m}R[\x ]}
\end{equation}
with $f_{i}\in F$ and $b_{i}x^{h_{i}}f_{i}\in R[\x ]_{\leq (-d,e)}$.
\end{definition}

An $F$-reducible polynomial belongs to $\bigcap _{m}
(P^{m}{R[\x] \! + \! I})$ but not necessarily to $I$.

\begin{proposition}\label{p:Grob-reduce}
A set $F$ satisfying \eqref{e:restriction} is a Gr\"obner basis of
$I= \langle F \rangle $ if and only if every element $p\in I$ is $F$-reducible.
\end{proposition}

\begin{proof}
Suppose every $p\in I$ is $F$-reducible.  Given $p\in I\setminus \{0
\}$, we are to show $\init (p)\in \langle \init (F) \rangle$.  Let
$\init (p) = ax^{e}$, $d = \ord (a)$.  If $x^{e}\not \in \langle \init
(F) \rangle$, then no summand in \eqref{e:reducible} has $e$ as the
exponent of its initial term, and hence every summand belongs to $R[\x
]_{<(-d,e)}$.  For $m>d$, this implies $p\in R[\x ]_{<(-d,e)}$, a
contradiction.

For the converse, fix an arbitrary $m$, and suppose $p\in I$ has no
expression of the form \eqref{e:reducible} for this $m$.  In
particular, $p\not \in P^{m}R[\x ]$, so $\init (p) = ax^{e}$, with
$\ord (a)<m$.  Since $d = \ord (a)$ is bounded above for all such $p$,
we may assume we have chosen $p$ so that $(-d,e)$ is minimal.  By
hypothesis there is some $f\in F$ such that $\init (f)$ divides
$x^{e}$, say $x^{e} = x^{h}\init (f)$.  Then $q = p-ax^{h}f$ has an
expression of the form \eqref{e:reducible} for this $m$, by the
minimality assumption.  But then so does $p$.
\end{proof}

\begin{remark}\label{r:graded-reduce}
Suppose $R[\x ]$ is given a multigrading $\deg \colon \N
^{n}\rightarrow A$, and $F$ consists of homogeneous polynomials.  Then
Proposition~\ref{p:Grob-reduce} holds in each degree separately: if
every nonzero $p\in I_{a}$ has $\init (p)\in \langle \init (F)
\rangle$, then every $p\in I_{a}$ is $F$-reducible.
\end{remark}

To each $f,g\in F$, there is an associated binomial syzygy $x^{u}\init
(f) = x^{v}\init (g) = \lcm (\init f,\init g)$.  We define the
corresponding {\em $S$-polynomial} as usual to be
\[
S(f,g) = x^{u}f - x^{v}g.
\]
Now we have a version of the Buchberger criterion for $F$ to be a
Gr\"obner basis.

\begin{proposition}\label{p:buch}
Let $B$ be a set of pairs $(f,g)\in F\times F$ whose associated
binomial syzygies generate the syzygy module of the initial terms
$\init (f)$, $f\in F$.  If $S(f,g)$ is $F$-reducible for all $(f,g)\in
B$, then $F$ is a Gr\"obner basis.
\end{proposition}

\begin{proof}
Fix $m\geq 0$.  We will show that every $p\in {I+P^{m}R[\x ]}$ has an
expression of the form \eqref{e:reducible} satisfying the conditions
in Definition~\ref{d:reducible} for this $m$.  We can assume $p\not
\in P^{m}R[\x ]$, so $\init (p)= ax^{e}$ with $d=\ord (a)<m$.
Certainly $p$ has {\it some} expression of the form
\eqref{e:reducible}, perhaps not satisfying $b_{i}x^{h_{i}}f_{i}\in
R[\x ]_{\leq (-d,e)}$.  Set $ x^{e_{i}} = \init (f_{i})$ and let
$(-d',e')$ be the maximum of $(-\ord (b_{i}),h_{i}+e_{i})$ over all
terms in our expression for $p$.  Since $p\not \in R[\x ]_{<(-d,e)}$,
we must have $(-d',e')\geq (-d,e)$.  In particular, $d'$ is bounded,
so we can assume our chosen expression for $p$ minimizes $(-d',e')$.
We are to show that $(-d',e') = (-d,e)$.

Suppose to the contrary that $(-d',e')>(-d,e)$.  Then we have $p\in
R[\x ]_{< (-d',e')}$, and every summand in \eqref{e:reducible} for which
$(-\ord (b_{i}), {h_{i}+e_{i}}) \not = (-d',e')$ is also in $R[\x ]_{<
(-d',e')}$.  Let $J$ be the set of indices $j$ for which $(-\ord
(b_{i}), {h_{i}+e_{i}}) = (-d',e')$.  The partial sum over these
indices in \eqref{e:reducible} must be in $R[\x ]_{<(-d',e')}$, so $\ord
(\sum _{J} b_{i})>d'$.  The Buchberger criterion implies that for all
indices $j,k\in J$, $x^{h_{j}}f_{j}-x^{h_{k}}f_{k}$ is a sum of
monomial multiples $x^{u}S(f,g)$ of $F$-reducible $S$-polynomials, all
satisfying $x^{u}\lcm (\init f,\init g) = x^{e'}$, and hence
$x^{u}S(f,g)\in R[\x ]_{\leq (0,e')}$.  Their $x^{e'}$ terms cancel, so in
fact they belong to $R[\x ]_{<(0,e')}$.  Being $F$-reducible, each
$x^{u}S(f,g)$ has an expression of the form \eqref{e:reducible} with
every term belonging to $R[\x ]_{< (0,e')}$, and hence so does
$x^{h_{j}}f_{j}-x^{h_{k}}f_{k}$.  Renaming indices so that $1\in J$,
we have
\[
\sum _{J} b_{i}x^{h_{i}}f_{i} = ({\textstyle \sum _{J}} b_{i})
x^{h_{1}} f_{1} + \sum _{J} b_{i}(x^{h_{i}}f_{i}-x^{h_{1}}f_{1}).
\]
The first term on the right belongs to $R[\x ]_{<(-d',e')}$, since
$h_{1}+e_{1} = e'$ and $\ord (\sum _{J}b_{i})>d'$.  In the second term
we can replace replace $b_{i}(x^{h_{i}}f_{i}-x^{h_{1}}f_{1})$ with an
expression of the form \eqref{e:reducible} with all terms in $R[\x ]_{<
(-d',e')}$.  Adding the remaining terms of our original expression for
$p$, we get a new expression with every term in $R[\x ]_{<(-d',e')}$.  This
contradicts the assumption that $(-d',e')$ was minimal.
\end{proof}

In order to apply the above results, we unfortunately need
the hypothesis \eqref{e:orderly}, which may fail in a non-Noetherian
ring.  We can still manage to avoid Noetherian hypotheses in
Theorem~\ref{t:multi-iso} by the device of reduction to the ground
ring $\Z $.  For this we need one last lemma, and then we will be
ready to prove our theorem.

\begin{lemma}\label{l:Hilb-over-Z}
Take $S$ and $h\colon A\rightarrow \N $ as in
Theorem~\ref{t:multi-Hilb}.  Then $H_{S_{D}}^{h}\cong (\Spec k)\times
_{\Z }H_{\Z [\x ]_{D}}^{h}$, for any subset $D$ of the degrees.  In
particular, $H_{S}^{h}\cong (\Spec k)\times _{\Z }H_{\Z [\x ]}^{h}$.
\end{lemma}

\begin{proof}
For simplicity, we only consider the case $D=A$, $S_{D}=S$.  The proof
in the general case is virtually identical.
Let $\hat{R}$ denote $R$ viewed as a $\Z $-algebra. Then
\begin{equation}\label{e:nat-bij}
H^{h}_{S}(R)= H_{\Z [\x ]}^{h}(\hat{R}),
\end{equation}
since we have the same set of ideals in $R[\x ] = {R\otimes _{k}S} =
{\hat{R}\otimes _{\Z }\Z [\x ]}$ on each side.  For any $X\in
\ZSchCateg $ and $R\in \kAlgCateg $, the $k$-morphisms $\Spec
R\rightarrow (\Spec k)\times _{\Z }X$ are in natural bijection with
the $\Z $-morphisms $\Spec R\rightarrow X$.  In other words, we have a
natural isomorphism $\fun{(\Spec k)\times X}(R)\cong
\fun{X}(\hat{R})$.  Together with \eqref{e:nat-bij}, this shows that
the $k$-schemes $H_{S}^{h}$ and $(\Spec k)\times _{\Z }H_{\Z [\x
]}^{h}$ have  isomorphic scheme functors.
\end{proof}

\begin{proof}[Proof of Theorem~\ref{t:multi-iso}] By
Lemma~\ref{l:Hilb-over-Z}, it is enough to prove the theorem in the
case $k=\Z$.  Since $D$ is supportive, the natural map
$H_{S}^{h}\rightarrow H_{S_{D}}^{h}$ is a closed embedding.  It
suffices to verify that it is an isomorphism locally on $
H_{S_{D}}^{h}$.  Specifically, let $U = \Spec R\subseteq
H_{S_{D}}^{h}$ be an affine open subset.  The closed embedding $U\cap
H_{S}^{h}\hookrightarrow U$ is given by a ring homomorphism
$R\rightarrow R/I$, and we are to show that $I=0$.  Localizing at
$P\in \Spec R$, it suffices to show that $I_{P}=0$.  The composite
morphism
\[
\Spec R_{P}\rightarrow \Spec R\rightarrow  H_{S_{D}}^{h}
\]
is an element of $H_{S_{D}}^{h} (R_{P})$.  We will show that it
belongs to the image of the map $H_{S}^{h} (R_{P})\hookrightarrow
H_{S_{D}}^{h} (R_{P})$.  This implies that the morphism $\Spec
R_{P}\rightarrow \Spec R$ factors through $\Spec R/I$, that is, the
localization homomorphism $R\rightarrow R_{P}$ factors through $R/I$
and therefore through $(R/I)_{P}$.  This yields a left inverse
$(R/I)_{P}\rightarrow R_{P}$ to the projection $R_{P}\rightarrow
(R/I)_{P}$, so $I_{P} = 0$.  Note that since we are assuming $k=\Z$,
and we have already shown that $H_{S_{D}}^{h}$ is quasiprojective
over $k$, the local ring $R_{P}$ is Noetherian, and hence satisfies
\eqref{e:orderly}.

We will show that the inclusion $H_{S}^{h} (R)\hookrightarrow
H_{S_{D}}^{h} (R)$ is surjective whenever $R$ is a local ring
satisfying \eqref{e:orderly}.  Let $L_{D}\subseteq R[\x ]_{D}$ be an
element of $H_{S_{D}}^{h} (R)$ and let $L'\subseteq R[\x ]$ be the
ideal generated by $L_{D}$.  By Remark~\ref{r:Fclosed}, we have
$L'_{D} = L_{D}$. In the proof of Theorem~\ref{t:multi-Hilb} we saw
that the conditions of Theorem~\ref{t:Hilb-scheme-II} hold.  We
conclude as in the proof of Theorem~\ref{t:Hilb-scheme-II} that $R[\x
]_{a}/L'_{a}$ is a finitely-generated $R$-module for all $a\in A$.
Let $K=R/P$ denote the residue field of $R$.  Then $KL'$ is the ideal
in $K[\x ]$ generated by $K{\otimes }L_{D}$.  Fix a term order on $\N
^{n}$ and let $J$ be the monomial ideal generated by the initial terms
$\{\init (p):p\in {K{\otimes }L_{D}} \}$.  For $a\in D$, $\dim K[\x
]_{a}/({K{\otimes }L_{a}}) = h(a)$.  Hence $J$ has Hilbert function
agreeing with $h$ on $D$, and by conditions (g) and (h) in
the definition of ``very supportive,'' $J$ has Hilbert function $h$.

The standard monomials for $J$ span $K[\x ]/KL'$.  By Nakayama's
Lemma, applied to each $R[\x ]_{a}/L'_{a}$ separately, it follows that
these standard monomials generate $R[\x ]/L'$ as an $R$-module.  What
remains to be shown is that they generate $R[\x ]/L'$ freely.  Then $L'$
is the required preimage of $L_{D}$ in $H_{S}^{h}(R)$.  For each
generator $x^{u}$ of the monomial ideal $J$, there is an element of
$KL'$ with initial term $x^{u}$.  Let $f\in L'$ be a representative of
this element modulo $PR[\x ]$.  The coefficient of $x^{u}$ in $f$ is a
unit in $R$, so we can assume it is $1$.  Then $\init (f) = x^{u}$.
Let $F$ be the set of polynomials $f$ obtained this way.

For $a\in D$, $R[\x ]_{a}/L'_{a}$ is free with basis the standard
monomials of degree $a$.  Given any monomial $x^{u}\in R[\x ]_{a}$,
its unique expansion modulo $L'_{a}$ by standard monomials belongs to
$R[\x ]_{\leq (0,u)}$.  To see this, observe that the expansion in
$K[\x ]$ of $x^{u}$ modulo $KL'$ contains only terms $x^{v}$ with
$v\leq u$.  It follows that the expansion of $bx^{u}$ belongs to $R[\x
]_{\leq (-\ord (b),u)}$.  Consider a nonzero element $p\in L'_{a}$,
with $\init (p)=bx^{e}$.  Replacing all remaining terms of $p$ with
their standard expansions, we get a polynomial $q\equiv p\pmod {L'}$.
At worst, this can change the coefficient of $x^{e}$ by an element of
$P^{\ord (b)+1}$, so $\init (q) = b'x^{e}$ for some $b'$.  All
remaining terms of $q$ are standard, and $q\in L'\setminus \{0 \}$, so
we must have $x^{e}\in J = \langle \init (F) \rangle$.  By
Remark~\ref{r:graded-reduce}, we deduce that every $p\in L'_{a}$ is
$F$-reducible.  In particular, $S(f,g)$ is $F$-reducible whenever the
generators $\init (f)$ and $\init (g)$ of $J$ participate in one of
the syzygies referred to in condition (s) for the very supportive set
$D$.  This shows that $F$ is a Gr\"obner basis for $I = \langle F
\rangle$.

Now, $I\subseteq L'$, and both $R[\x ]_{D}/I_{D}$ and $R[\x
]_{D}/L'_{D}$ are free with basis the standard monomials in degrees
$D$, so $I_{D} = L'_{D}$.  Both $I$ and $L'$ are generated in degrees
$D$, so $I = L'$.  Hence the standard monomials are $R$-linearly
independent modulo $L'$.
\end{proof}

When the grading is positive and the Hilbert scheme is projective, the
preceding results lead to an explicit description of the multigraded
Hilbert scheme $H_{S}^{h}$ by equations in Pl\"ucker coordinates,
although the number of variables and equations involved may be
extremely large.  We write $a\leq b$ for degrees $a,b\in A$ if $b-a\in
A_{+}$.  Since our grading is positive, this is a partial ordering on
the degrees.  For any finite set of degrees $D \subseteq A$, the
Hilbert functor $H_{S_D}^{h}$ is defined as a subfunctor of the
Grassmann functor $G^h_{S_D}$ by the conditions on $L\in
H_{S_{D}}^{h}(R)$:
\begin{equation} \label{quadraticEquations}
\text{for all $a<b \in D$ and all $x^{u}$ with $\deg (u)=b-a$: $x^{u}
L_a \subseteq L_b$.}
\end{equation}
For a positive grading, there are finitely many monomials in each
degree.  Each member of the above finite system of inclusions
translates into well-known quadratic equations in terms of Pl\"ucker
coordinates on $G^{h(a)}_{S_{a}}\times G^{h(b)}_{S_{b}}$.  Together
these equations describe the Hilbert scheme $H_{S_{D}}^{h}$ as a
closed subscheme of $G^h_{S_D}$.  We call \eqref{quadraticEquations}
the {\em natural quadratic equations}.

\begin{corollary}\label{c:quadrics}
If the grading is positive and  $D \subseteq A$ is
very supportive for $h \colon A \rightarrow \N$ then the Hilbert
scheme $H_{S}^{h}$ is defined by the natural quadratic equations
\eqref{quadraticEquations}.
\end{corollary}

Let $D \subseteq E$ be two finite sets of degrees, where $D$ is
supportive and $E$ is very supportive. Then our problem is to write
down equations for the image of the closed embedding of $H^h_S \simeq
H^h_{S_E} $ into $H^h_{S_D}$ given by Corollary \ref{p:multi-closed}.
Each degree $e\in E \backslash D$ contributes to these equations,
which we have already described in the discussion following the proof
of Theorem~\ref{t:Hilb-scheme-II} as the Fitting ideal for a certain
matrix $\Gamma $.  In the positively graded case, this matrix is
finite and we can describe it explicitly.  The columns of $\Gamma $
correspond to the monomials of degree $e$.  For each degree $d\in D$,
$d<e$, and each set $B$ consisting of $h(d)+1$ monomials of degree
$d$, there is an element $\sum _{b\in B}{\gamma _{B\setminus \{b
\}}\cdot b}$ of $L_{d}$, where $\gamma _{B\setminus \{b \}}$ denotes
the Pl\"ucker coordinate on $G^{h(d)}_{S_{d}}$indexed by the set of
$h(d)$ monomials $B\setminus \{b \}$.  Equation \eqref{e:bayer-coords}
in Example~\ref{daveRunning} illustrates this.  Multiply each such
generator of $L_{d}$ by a monomial $x^{u}$ of degree $e-d$ to get a
homogeneous polynomial of degree $e$ in $\x $ with coefficients that
are Pl\"ucker coordinates.  The vector of coefficients gives a row of
$\Gamma $, which is the matrix of all rows obtained in this way.
Setting $r = \rank S_{e} = \binom{n+e-1}{e}$, the minors
\begin{equation}\label{e:determinantal}
I_{r-h(e)+1}(\Gamma )
\end{equation}
are the {\em natural determinantal equations} contributed by the
degree $e$.

\begin{theorem}\label{c:determinants}
If $D \subseteq A$ is supportive for $h \colon A \rightarrow \N$, then
the Hilbert scheme $H_{S}^{h}$ is defined by the natural quadratic
equations \eqref{quadraticEquations} and the natural determinantal
equations \eqref{e:determinantal}, where $e$ runs over 
$E\setminus D$, for a very supportive superset $E$ of $D$.
\end{theorem}

\section{The Grothendieck Hilbert scheme}
\label{sec4}

In this section we relate our construction to
Grothendieck's classical Hilbert scheme.  Expressing the latter as a
special case of the multigraded Hilbert scheme, 
our natural quadratic equations will become 
Gotzmann's equations \cite{Got}, while the natural determinantal
equations become those of Iarrobino and Kleiman \cite{IK}.  Dave Bayer
in his thesis \cite[\S VI.1]{Bay} proposed a more compact
system of determinantal equations, each having degree $n$ in 
Pl\"ucker coordinates, and he conjectured that they also
define Grothendieck's Hilbert scheme.  
Here we prove Bayer's conjecture.

The {\em Grothendieck Hilbert scheme} $\Hilb _{n-1}^{g}$ represents
the functor of flat families $X\subseteq \P ^{n-1}(R)$, $R\in
\kAlgCateg $, with a specified Hilbert polynomial $g$.   The
homogeneous coordinate ring of $\P ^{n-1}(R)$ is $R[\x ] =
R[x_{1},\ldots,x_{n}]$, and the ideal of $X$ is a saturated
homogeneous ideal $L\subseteq R[\x ]$ such that in sufficiently large
degrees, $R[\x ]/L$ is locally free with Hilbert function $g$.  
Let $d_0 = d_0(g,n)$ denote the {\it Gotzmann number} 
\cite[Definition C.12]{IK}. Gotzmann \cite{Got} proved:
(1) every saturated ideal with Hilbert polynomial $g$ has Hilbert
function  $g$ in degrees $d\geq d_{0}$, and (2) every
ideal with Hilbert function $g$ in degrees $d\geq d_{0}$ coincides in
these degrees with its saturation.  

\begin{lemma} Grothendieck's Hilbert scheme $\Hilb _{n-1}^{g}$
is isomorphic to the multigraded Hilbert scheme
 $H _{S}^{h}$, where $S = k[\x ]$ with the standard
$\Z$-grading, with Hilbert function $h$ defined by
 $h(d) = g(d)$ for $d\geq d_{0}$, $h(d) =
\binom{n+d-1}{d}$ for $d<d_{0}$. 
\end{lemma}

\begin{proof}
The ideals described by the functor $H_{S}^{h}$ are the
truncations to degrees $d\geq d_{0}$ of the ideals described by the
Grothendieck functor $\Hilb_{n-1}^{g}$.  A natural bijection between
the two is given by truncation in one direction and saturation in the
other. Hence both schemes represent the same functor.
\end{proof}

The Gotzmann number $d_{0}$ equals the maximum of the
Castelnuovo-Mumford regularity of any saturated monomial ideal $I$
with Hilbert polynomial $g$ \cite[Proposition C.24]{IK}.  The set of
such ideals is finite by Lemma~\ref{l:Diane}.  For a monomial ideal,
the regularity of $I$ is a purely combinatorial invariant, equal to
the maximum over all $i$ and all minimal $i$-th syzygies of $d-i$,
where $d$ is the degree of the syzygy.  The regularity will not exceed
$d_0$ if $I$ is replaced by its truncation to degrees $\geq d_{0}$.
It follows that for every monomial ideal $I$ generated in degree
$d_{0}$ and with Hilbert polynomial $g$, the Hilbert function of $I$
coincides with $g$ in degrees $\geq d_{0}$, and $I$ has a linear free
resolution.  In particular the minimal $S$-pairs of $I$ have degree
$d_{0}+1$.  These considerations show that Gotzmann's Regularity
Theorem and Persistence Theorem can be rephrased in the language of
the previous section as follows:

\begin{proposition} \label{GotzPersis}
Let $g$ be a Hilbert polynomial, defining $d_0$ and $h$ as above.
Then $D =\{d_0\}$ is supportive and $E = \{d_0,d_0+1\}$ is very
supportive for $H _{S}^{h} = \Hilb _{n-1}^{g}$.
\end{proposition}

We can now write equations for the Grothendieck Hilbert scheme in two
possible ways. The set $E = \{d_0,d_0+1\}$ gives an embedding into a
product of Grassmannians
\begin{equation} \label{TwoGrass}
H_S^{h} = H_{S_E}^{h} \hookrightarrow G^{h(d_0)}_{S_{d_0}}
\times G^{h(d_0+1)}_{S_{d_0+1}} .
\end{equation}
This is the embedding described by Gotzmann in \cite[Bemerkung
(3.2)]{Got}; see also \cite[Theorem C.29]{IK}. It is defined
scheme-theoretically by the {\em natural quadratic equations} given in
\eqref{quadraticEquations}.  We illustrate these equations with a
simple example.

\begin{example} \label{TwoPointsInThePlane}
Take $S = k[x,y,z]$ with Hilbert function $h(0)=1$ and $h(d) = 2 $ for
$d \geq 1 $.  Our Hilbert scheme $ H_S^{h} $ coincides with the
Grothendieck Hilbert scheme of two points in the projective plane
$\P^2$. The Gotzmann number is $d_0 = 2$. The pair $E = \{2,3\}$ is
very supportive and gives the embedding \eqref{TwoGrass} into the
product of Grassmannians $G^2_6 \times G^2_{10}$.  The Pl\"ucker
coordinates for the Grassmannian $G^2_6$ (resp.~$G^2_{10}$) are
ordered pairs of quadratic (resp.~cubic) monomials in $x,y,z$. These
define the Pl\"ucker embeddings $G^2_6 \hookrightarrow \P^{14} $ and
$G^2_{10} \hookrightarrow \P^{44} $.  The Hilbert scheme $H_S^{h}$ is
the closed subscheme of $G^2_6 \times G^2_{10}$ defined by $600$
bilinear equations as in \eqref{quadraticEquations}.  There are $180$
two-term relations such as
\[
[x y^2, x y z] \cdot [y z, x y] + [x y^2, y^2 z] \cdot [x y, x z]
= 0,
\]
and $420$ three-term relations such as
\[
[x^2 z, x y^2] \cdot [x z, y z] + [x^2 z, x y z] \cdot [y z, x y] +
[x^2 z, y^2 z] \cdot [x y, x z] = 0.
\]
The validity of these equations is easily checked for subschemes of
$\P^2$ consisting of two distinct reduced points $(x_1:y_1:z_1) $ and
$(x_2:y_2:z_2) $.  Just replace each bracket by the corresponding $2
\times 2$ determinant, as in $ [x^2 z, x y^2] \mapsto x_1^2 z_1 x_2
y_2^2 - x_2^2 z_2 x_1 y_1^2 $. \qed
\end{example}

\smallskip

In the remainder of this section we will study not the Gotzmann
embedding \eqref{TwoGrass} but the other (more efficient) embedding
given by Proposition \ref{GotzPersis}.  Namely, the supportive
singleton $D = \{d_0\}$ defines the embedding into a single
Grassmannian
\begin{equation}
\label{OneGrass} H_S^{h} \hookrightarrow H_{S_D}^{h} =
G^{h(d_0)}_{S_{d_0}}.
\end{equation}
This is the embedding described in Bayer's thesis \cite[\S VI.1]{Bay}
and in \cite[Prop.~C.28]{IK}.  It follows from Theorem
\ref{c:determinants} that the Hilbert scheme is defined as a closed
subscheme of the Grassmannian by the {\em natural determinantal
equations} \eqref{e:determinantal}.  Iarrobino and Kleiman proved this
in the present case in \cite[Proposition C.30]{IK}, so we refer to the
equations \eqref{e:determinantal} for the Grothendieck Hilbert scheme
as the {\it Iarrobino--Kleiman equations}.  Note that the
Iarrobino--Kleiman equations for the embedding \eqref{OneGrass} are
homogeneous polynomials of degree $\binom{n+d_{0}}{d_{0}+1} -
h(d_{0}+1) + 1$ in the Pl\"ucker coordinates.

We now present a third system of homogeneous equations for the
Grothendieck Hilbert scheme, which Bayer proved define it
set-theoretically.  Like the Iarrobino--Kleiman equations, Bayer's
equations are homogeneous equations in the Pl\"ucker coordinates on
the single Grassmannian $G_{S_{d_{0}}}^{h(d_{0})}$.  However, Bayer's
equations are more compact: their degree always equals $n$, the number
of variables, independently of $g$, $h$ and $d_{0}$
\cite[p.~144]{Bay}.  Bayer conjectured that his equations define the
correct scheme structure \cite[p.~134]{Bay}.  We will prove this
conjecture.

\begin{theorem} \label{BayerConj}
Grothendieck's Hilbert scheme parametrizing subschemes of $\P ^{n-1} $
with any fixed Hilbert polynomial is defined in the Grassmannian
embedding \eqref{OneGrass} by Bayer's equations, which are homogeneous
of degree $n$ in the Pl\"ucker coordinates.
\end{theorem}

Although the Bayer equations define the same subscheme of the
Grassmannian as the Iarrobino--Kleiman equations, they do not generate
the same homogeneous ideal.  This phenomenon is hardly surprising,
since any projective scheme can be defined by many different
homogeneous ideals.  Even the Bayer equations are often not the
simplest ones: the common saturation of both ideals frequently
contains equations of degree less than $n$.  This happens for Example
\ref{TwoPointsInThePlane}, which will be reexamined below, and it
happens for \cite[Example C.31]{IK}, where the Iarrobino--Kleiman
equations have degree $25$ while the Bayer equations have degree $3$.

The best way to introduce Bayer's equations and relate them to the
Iarrobino--Kleiman equations is with the help of Stiefel coordinates
on the Grassmannian.  For the remainder of this section we use the
following abbreviations:
\[
 d = d_0;\; h = h(d);\; h' = h(d+1);\; r = \binom{n+d-1}{d};\; r' =
\binom{n+d}{d+1}.
\]
As before, $G^{h}_{r}$ denotes the Grassmann scheme parametrizing
quotients of rank $h$ of $S_d$.  We digress briefly to review the
relationship between local coordinates, Stiefel coordinates, and
Pl\"ucker coordinates.

Recall from Section~\ref{sec2} that the Grassmannian $G^{h}_{r}$ is
covered by affine charts $G^{h}_{r\setminus B}$, whose functor
$\fun{G^{h}_{r\setminus B}}(R)$ describes free quotients $R^{r}/L$
with basis $B$, where $B$ is an $h$-element subset of some fixed basis
$X$ of $k^{r}$.  Here we identify $k^{r}$ with $S_{d}$, and $X$ is the
set of all monomials of degree $d$.  At a point $L\in
\fun{G^{h}_{r\setminus B}}(R)$, the local (affine) coordinates $\gamma
^{x}_{b}$ take unique values in $R$ such that
\[
x\; \equiv \; \sum _{b\in B}\gamma ^{x}_{b}\cdot b {\pmod L} \quad
\text{for all $x\in X\setminus B$}.
\]
Consider the $(r-h)\times r$ matrix $\Gamma $ with columns indexed by
the elements of $X$, constructed as follows.  Index the rows of
$\Gamma $ by the elements of $X \backslash B$.  In the column indexed by
$b\in B$, put the coordinates $-\gamma ^{x}_{b}$ of $L$, for $x\in
X \backslash B$.  In the complementary square submatrix with columns
indexed by $X\setminus B$, put an $(r-h)\times (r-h)$ identity matrix.
Then the rows of $\Gamma $ span the submodule $L\subseteq R^{r}$.

More invariantly, if we insist that $R^{r}/L$ be free, not just
locally free, but do not choose the basis $B$ in advance, we can
always realize $L$ as the row space of some $(r-h)\times r$ matrix
$\Omega $, at least one of whose maximal minors is invertible in $R$.
The entries of $\Omega $ are the {\it Stiefel coordinates} of $L$.
They are well-defined up to change of basis in $L$, that is, up to
multiplication of $\Omega $ on the left by matrices in $GL_{r-h}(R)$.
A little more generally, we can regard any $(r-h)\times r$ matrix
$\Omega $ whose maximal minors generate the unit ideal in $R$ as the
matrix of Stiefel coordinates for its row-space $L\subseteq R^{r}$, as
$R^{r}/L$ will then be locally free of rank $h$.

When $R^{r}/L$ is locally free of rank $h$, its top exterior power
$\wedge ^{h}(R^{r}/L)$ is a rank-$1$ locally free quotient of $\wedge
^{h} (R^{r})$, corresponding to an element of $\fun{\P
^{\binom{r}{h}-1}}(R)$.  The Pl\"ucker embedding
$G^{h}_{r}\hookrightarrow \P ^{\binom{r}{h}-1}$ is given in scheme
functor terms by the natural transformation sending $L$ to the kernel
of $\wedge ^{h} (R^{r}) \rightarrow \wedge ^{h}(R^{r}/L)$.  The
homogeneous coordinates on $\P ^{\binom{r}{h}-1}$ are {\it Pl\"ucker
coordinates}.  They are indexed by exterior products of the elements
of $X$ and denoted
\begin{equation}\label{Brackets}
[x_{1},\ldots,x_{h}].
\end{equation}
In terms of Stiefel coordinates, we can identify
$[x_{1},\ldots,x_{h}]$ with the maximal minor of $\Omega $ whose
columns are indexed by $x_{1},\ldots,x_{h}$, up to a sign depending on
the order of the monomials in the bracket.

Some caution is due when using Stiefel and Pl\"ucker coordinates in
the scheme functor setting: for an arbitrary $L\in
\fun{G^{h}_{r\setminus B}}(R)$, the matrix $\Omega $ of Stiefel
coordinates need not exist, as $L$ may not be generated by $r-h$
elements.  This difficulty arises even for homogeneous coordinates on
projective space (the special case $h=r-1$).  Nevertheless, for the
purpose of determining the ideal of a closed subscheme $H\subseteq
G^{h}_{r}$, it suffices to consider the restriction of the scheme
functors involved to {\it local} rings $R$.  Stiefel and Pl\"ucker
coordinates then make sense for any $R$-valued point $L$.  Throughout
the rest of this section, $R$ will always denote a local ring.

The basic observation leading to the Bayer equations is that when a
subscheme of $G^{h}_{r}$ is defined by nice enough equations in
Stiefel coordinates, they can sometimes be converted to equations of
much lower degree in Pl\"ucker coordinates. For instance, the
submodule $L_2\subseteq R[x,y,z]_{2}$ in
Example~\ref{TwoPointsInThePlane} is spanned by four quadrics,
\begin{eqnarray*}
& { a_1}{x}^{2}+{ a_2}xy+{ a_3}xz+{ a_4}{y}^{2}+
{ a_5}yz+{ a_6}{z}^{2} ,\\
& { b_1}{x}^{2}+{ b_2}xy+{
 b_3}xz+{ b_4}{y}^{2}+{ b_5}yz+{ b_6}{z}^{2}, \\
& { c_1}{x}^{2}+{ c_2}xy+{ c_3}xz+{ c_4}{y}^{2}+{
 c_5}yz+{ c_6}{z}^{2}, \\
& { d_1}{x}^{2}+{ d_2}xy+{
 d_3}xz+{ d_4}{y}^{2}+{ d_5}yz+{ d_6}{z}^{2}.
\end{eqnarray*}
The matrix $\Omega$ is the $4 \times 6$ matrix of coefficients, which
are the Stiefel coordinates.  The fifteen $4 \times 4$ minors of
$\Omega$ are identified with the fifteen Pl\"ucker coordinates on
$G^2_6$. Some care is required with the signs; for
instance,
\begin{equation}  \label{fromTwotoFour}
[yz, z^2] = \det \begin{pmatrix}
 a_{1}&	a_{2} &	a_{3}&	a_{4}\\
 b_{1}& b_{2}&	b_{3} & b_{4}\\
 c_{1}& c_{2}&	c_{3} & c_{4}\\
 d_{1}& d_{2}&	d_{3} & d_{4}\\
\end{pmatrix},
\quad [y^2, z^2] = - \det \begin{pmatrix}
 a_{1}&	a_{2} &	a_{3}&	a_{5}\\
 b_{1}& b_{2}&	b_{3} & b_{5}\\
 c_{1}& c_{2}&	c_{3} & c_{5}\\
 d_{1}& d_{2}&	d_{3} & d_{5}\\
\end{pmatrix}.
\end{equation}

Returning to the general discussion, observe that the image $\x \cdot
L_d$ of $L_d\subseteq R[\x ]_{d}$ is spanned by $ x_1 L_d , x_2 L_d,
\ldots, x_n L_d $ inside $R[\x ]_{d+1} = R^{r'}$.  We may represent
$\x \cdot L_d$ by a matrix $\widehat \Omega $ with $n(r-h)$ rows and
$r'$ columns.  The matrix $\widehat \Omega$ contains $n$ copies of the
matrix $\Omega$ and is otherwise zero. The columns of $\widehat
\Omega$ are labeled by the monomials in $R[\x ]_{d+1}$ in
lexicographic order. In our running example, we have
\begin{equation} \label{LeftMatrix}
\widehat\Omega = 
\left (\begin {array}{cccccccccc}
{ a_1}&{ a_2}&{ a_3}&{ a_4}&{ a_5}&{ a_6}&0&0&0&0\\
{ b_1}&{b_2}&{ b_3}&{ b_4}&{ b_5}&{ b_6}&0&0&0&0 \\
{ c_1}&{ c_2}&{ c_3}&{ c_4}&{c_5}&{ c_6}&0&0&0&0\\
{ d_1}&{ d_2}&{d_3}&{ d_4}&{ d_5}&{ d_6}&0&0&0&0\\
0&{a_1}&0&{ a_2}&{ a_3}&0&{ a_4}&{ a_5}&{ a_6}&0 \\
0&{ b_1}&0&{ b_2}&{ b_3}&0&{ b_4}&{b_5}&{ b_6}&0\\
0&{ c_1}&0&{ c_2}&{c_3}&0&{ c_4}&{ c_5}&{ c_6}&0\\
0&{d_1}&0&{ d_2}&{ d_3}&0&{ d_4}&{ d_5}&{ d_6}&0 \\
0&0&{ a_1}&0&{ a_2}&{ a_3}&0&{ a_4} &{ a_5}&{ a_6}\\
0&0&{ b_1}&0&{ b_2}&{b_3}&0&{ b_4}&{ b_5}&{ b_6}\\
0&0&{c_1}&0&{ c_2}&{ c_3}&0&{ c_4}&{ c_5}&{ c_6} \\
0&0&{ d_1}&0&{ d_2}&{ d_3}&0&{ d_4} &{ d_5}&{ d_6}
\end {array}\right )
\end{equation}
with columns labeled
$x^3, x^2 y, x^2 z, x y^2, x y z, x z^2, y^3, y^2 z, y z^2, z^3$.

The choice of $d$ as the Gotzmann number ensures that $\widehat{\Omega
}$ has an invertible minor of order $r'-h'$ whenever $\Omega $ has an
invertible maximal minor.  The natural determinantal equations
\eqref{e:determinantal} defining $H^{h}_{S}$ as a closed subscheme of
$ G^{h}_{r} $ are the minors of order $r'-h'+1$ of the matrix
$\widehat \Omega$.  They are the Iarrobino--Kleiman equations
expressed in Stiefel coordinates, and are exactly the equations which
ensure that $R[\x ]_{d+1}/\x \cdot L_{d}$ is locally free of rank
$h'$.  In our example, we are looking at $2,200 = \binom{12}{9} \times
\binom{10}{9}$ polynomials of degree $9$.  We wish to replace these by
a smaller number of cubic polynomials in the $4 \times 4$ minors of
the matrix $\Omega$.

In general, our problem is this: Let $J$ be the Fitting ideal
generated by the minors of order $r'-h' + 1$ of the matrix $\widehat
\Omega$.  This is an ideal in the polynomial ring $k[\Omega]$
generated by entries of $\Omega$, that is, by the Stiefel coordinates,
viewed as indeterminates.  We seek an ideal $J'$ generated by
polynomials of degree $n$ in the Pl\"ucker coordinates, or maximal
minors of $\Omega $, such that $J$ and $J'$ define systems of
equations which have the same solutions $\Omega $ over any local ring
$R$.

We now give Bayer's construction and show that it solves the above
problem.  Let ${\Omega \otimes S_1}$ be the matrix representing the
submodule ${S_1 \otimes_{k} L_d}$ of the tensor product ${S_1
\otimes_{k} R[\x ]_{d}}$. Thus ${\Omega \otimes S_1}$ is a matrix with
$n (r-h)$ rows and $n r$ columns.  The row labels of ${\Omega \otimes
S_1}$ coincide with the row labels of $\widehat \Omega$. We form their
concatenation
$$
\bigl( \widehat \Omega \mid \Omega \otimes S_1 \bigr).
$$
Bayer's equations are certain maximal minors of this matrix.  Each
column of $ \widehat \Omega $ is a sum of columns of $\Omega \otimes
S_1 $, and these sums involve distinct leading columns.  Therefore we
may---for the sake of efficiency---pick a submatrix $(\Omega \otimes
S_1)_{\red} $ of $(\Omega \otimes S_1) $ of format $n (r - h) \times (n
r -r')$ such that the maximal minors of
\begin{equation}
\label{DavesMatrix} \bigl( \widehat \Omega \mid ( \Omega \otimes
S_1)_{\red} \bigr)
\end{equation}
have the same $\Z$-linear span as those of $\bigl( \widehat \Omega
\mid \Omega \otimes S_1 \bigr)$.  Note that the matrix
\eqref{DavesMatrix} has $n (r - h)$ rows and $n r$ columns.  Each
maximal minor of \eqref{DavesMatrix} is a homogeneous polynomial of
degree $n(r - h)$ in $k[\Omega]$, and, by Laplace expansion, it can be
written as a homogeneous polynomial of degree $n $ in the Pl\"ucker
coordinates \eqref{Brackets}.  The {\it Bayer equations} are those
maximal minors of \eqref{DavesMatrix} gotten by taking any set of
$r'-h' + 1$ columns of $ \widehat \Omega$ and any set of $n (h-r) -
r'+h' - 1$ columns of $({\Omega \otimes S_1})_{\red} $.

In our running example, we take the reduced tensor product matrix
as follows:
\begin{equation} 
\label{RightMatrix}
( \Omega \otimes S_1)_{\red} = 
\left (\begin {array}{cccccccc}
{ a_2}&{ a_3}&{ a_4}&{a_5}&0&{ a_6}&0&0\\
{ b_2}&{ b_3}&{b_4}&{ b_5}&0&{ b_6}&0&0\\
{ c_2}&{c_3}&{ c_4}&{ c_5}&0&{ c_6}&0&0\\
{d_2}&{ d_3}&{ d_4}&{ d_5}&0&{ d_6}&0&0 \\
0&0&0&0&{ a_3}&0&{ a_5}&{ a_6} \\
0&0&0&0&{ b_3}&0&{ b_5}&{ b_6} \\
0&0&0&0&{ c_3}&0&{ c_5}&{ c_6} \\
0&0&0&0&{ d_3}&0&{ d_5}&{ d_6} \\
0&0&0&0&0&0&0&0\\
0&0&0&0&0&0&0&0 \\
0&0&0&0&0&0&0&0\\
0&0&0&0&0&0&0&0
\end {array}\right )
\end{equation}
The matrix (\ref{DavesMatrix}) has format $12 \times 18$, and each of
its maximal minors is a homogeneous polynomial of degree $3$ in the
$15$ Pl\"ucker coordinates $[x^2, xy]$, $[x^2 ,xz]$, $\ldots$, $[yz,
z^2]$.  There are $560 = \binom{10}{9} \times \binom{8}{3}$ Bayer
equations, gotten by taking any $9$ columns from \eqref{LeftMatrix}
and any $3$ columns from \eqref{RightMatrix}.

\begin{proof}[Proof of Theorem \ref{BayerConj}:] Clearly, every Bayer
equation belongs to the Fitting ideal $I_{r'-h'+1}(\widehat{\Omega
})$.  We must show that (for $R$ local) the vanishing of the Bayer
minors implies that $I_{r'-h'+1}(\widehat{\Omega })=0$.  This would be
obvious if the matrix $({\Omega \otimes S_{1}})_{\red }$ contained an
identity matrix as a maximal square submatrix.  But the Bayer ideal is
unchanged if we use ${\Omega \otimes S_{1}}$ in place of $({\Omega
\otimes S_{1}})_{\red }$, and it is $GL_{n(r-h)}(R)$-invariant.  Hence
it suffices that ${\Omega \otimes S_{1}}$ have some maximal minor
invertible in $R$.  This follows from the fact that $\Omega $ has such
a minor.
\end{proof}

While the Bayer equations do define the correct scheme structure on
the Hilbert scheme, they are far from minimal with this property.  For
instance, in our example, there are $560$ Bayer cubics which, together
with the $15$ quadratic Pl\"ucker relations for $G^2_6$, define the
Hilbert scheme $H^{h}_{S}$ as a closed subscheme of dimension $4$ and
degree $21$ in $\P^{14}$.  However, $H^{h}_{S}$ is irreducible and its
prime ideal is the ideal of algebraic relations on the $2 \times 2$
minors of the matrix
$$
\left (\begin {array}{cccccc}
x_1^2 &x_1 y_1 &x_1 z_1 & y_1^2 &y_1 z_1 & z_1^2  \\
x_2^2 &x_2 y_2 &x_2 z_2 & y_2^2 &y_2 z_2 & z_2^2  
\end {array} \right)
$$
This prime ideal is minimally generated by $45$ quadrics.

\section{Toric Hilbert schemes and their Chow morphisms }
\label{sec5}

In this section we examine Hilbert schemes which arise in toric
geometry.  Our goals are to describe equations for the {\em toric
Hilbert scheme}, and to define the {\em toric Chow morphism}.  In the
process we answer some questions left open in earlier investigations
by Peeva, Stillman, and the second author.  We fix an $A$-grading of
the polynomial ring $S = k[\x ]$ and consider the constant Hilbert
function
\begin{equation} \label{constantFct}
h(a) = 1 \quad \text{for all $a\in A_{+}$}.
\end{equation}
The multigraded Hilbert scheme $H^{1}_{S}$ defined by this Hilbert
function is called the {\em toric Hilbert scheme}.  Its functor
$H^{1}_S(R)$ parametrizes ideals $I\subseteq R[\x ]$ such that $(R[\x
]/I)_a$ is a rank-one locally free $R$-module for all $a \in A_+$.

Assuming that the elements $a_{i} = \deg (x_{i})$ generate $A$, we
have a presentation
\begin{equation}\label{e:Apresentation}
0 \rightarrow M \rightarrow \Z ^{n} \rightarrow A \rightarrow 0,
\end{equation}
which induces a surjective homomorphism of group algebras over $k$,
\begin{equation}\label{e:homomorphismA->T}
k[\x ,\x ^{-1}] = k [ \Z ^{n}] \rightarrow k [ A],
\end{equation}
and a corresponding closed embedding of $G = \Spec k [A] $ as an
algebraic subgroup of the torus $\T ^{n} = \Spec k[\x ,\x ^{-1}]$.
The torus $\T ^{n}$ acts naturally on $\A ^{n}$ as the group of
invertible diagonal matrices, and so its subgroup $G$ also acts on $\A
^{n}$.  An ideal $I\subseteq R\otimes S$ is homogeneous for our
grading if and only if the closed subscheme defined by $I$ in $\A
^{n}_{R} = \A ^{n}\times _{k}\Spec R$ is invariant under the action of
$G_{R} = G\times _{k}\Spec R$.

\begin{remark} If $A$ is a finite abelian group then the toric Hilbert
scheme $H^{1}_S$ coincides with Hilbert scheme $\Hilb ^{G}(\A ^{n})$
of regular $G$-orbits studied by Nakamura \cite{N}.  If the group $A$
is free abelian and the grading is positive, then $H^{1}_S$ coincides
with the toric Hilbert scheme studied by Peeva and Stillman
\cite{PSlocal, PS}. \qed
\end{remark}

There is a distinguished point on the toric Hilbert scheme $H^{1}_S$,
namely, the ideal
\begin{equation*}
I_M \quad = \quad \langle \, x^u - x^v  \,\, : \,\,\,
u , v \in \N^n , \,
  \deg (u)=\deg (v) \, \rangle.
\end{equation*}
Note that $\deg (u)=\deg (v)$ means that $u-v$ lies in the sublattice
$M$ in \eqref{e:Apresentation}.  Restricting the ring
map in \eqref{e:homomorphismA->T} to $S = k[\N ^{n}]$, its
kernel is $I_{M}$.  Hence, identifying $\A ^{n} = \Spec S$ with the
space of $n\times n$ diagonal matrices, and $\T ^{n}$ with its open
subset of invertible matrices, $I_{M}$ is the ideal of the closure in
$\A ^{n}$ of the subgroup $G\subseteq \T ^{n}$.

A nonzero binomial $x^u - x^v \in I_{M}$ is called  {\em Graver} if there is no
other binomial $x^{u'} - x^{v'} $ in $I_M$ such that $x^{u'}$ divides
$x^{u}$ and $x^{v'}$ divides $x^{v}$.  The degree $a = \deg(u) =
\deg(v)$ of a Graver binomial is a {\em Graver degree}.  The
set of Graver binomials is finite. The finite set of  all Graver
degrees can be computed using Algorithm 7.2 in \cite{StuBook}.

\begin{proposition}\label{p:graver}
The set of Graver degrees is supportive, and the natural determinantal
equations \eqref{e:determinantal} for this set coincide with the
determinantal equations for the toric Hilbert scheme given by Peeva
and Stillman in \cite[Definition 3.3]{PS}.
\end{proposition}

\begin{proof} 
The Graver degrees are supportive by \cite[Proposition 5.1]{PSlocal};
the proof given there for positive gradings works for nonpositive
gradings as well.  The Fitting equations in \cite[Definition 3.3]{PS}
are precisely our Fitting equations \eqref{e:determinantal}, in the
special case when the Hilbert function $h$ is the constant $1$.
\end{proof}

In the positively graded case, a doubly-exponential bound was given in
\cite[Proposition 5.1]{S} for a set of degrees which is very
supportive for the toric Hilbert scheme.  Peeva \cite[Corollary
5.3]{PSlocal} improved the bound to single-exponential and gave an
explicit description of a very supportive set $E$ in \cite[Theorem
5.2]{PSlocal}.

\begin{proposition}\label{p:quad-binoms}
Let $\deg \colon \N ^{n}\rightarrow A$ be a positive grading and $E
\subseteq A_+$ a finite, very supportive set of degrees for the toric
Hilbert scheme. Then the natural quadratic equations
\eqref{quadraticEquations} are precisely the quadratic binomials given
in \cite[Equation (5.3)]{S}.
\end{proposition}

\begin{proof} We only need to make explicit the equations expressing
condition \eqref{quadraticEquations}:
\[
\text{for all $a<b \in E$ and all $x^{u}$ with $\deg (u)=b-a$: $x^{u}
L_a \subseteq L_b$.}
\]
Let $R$ be a local ring.  For the constant Hilbert function $h=1$, the
ambient graded Grassmann scheme is a product of projective spaces, one
for each degree:
\[
G^{1}_{S_{E}} = \prod_{a \in E} \P( S_a).
\]
For each monomial $x^u$ in $S_a$ there is a coordinate $z^a_u$ on the
projective space $\P(S_a)$, such that the
$z^a_u$ for $\deg (u) = a$ are the Pl\"ucker coordinates on $\fun{\P
(S_{a})}(R)$.  The submodule $L_a$ of ${R\otimes S_a}$ represented by
a point $(z^a_u)$ in $\fun{\P(S_a)}(R)$ is generated by
\[
z^a_u \cdot x^v - z^a_v \cdot x^u \quad \text{for all $\deg (u)=\deg
(v)=a$}.
\]
For $R$ local, condition \eqref{quadraticEquations} is thus equivalent
to the system of binomial equations
\begin{equation} 
\label{BerndsOldBinomials} z^a_u \cdot z^b_{v+w} = z^a_v \cdot
z^b_{u+w} \quad \text{for $a,b\in E$, $\deg (u)=\deg (v)=a$, $\deg (w)
= b-a$},
\end{equation} 
which are precisely the equations in \cite[(5.3)]{S}.  A closed
subscheme cut out by equations in any scheme is determined by the
evaluation of its subfunctor on local rings $R$.  Hence $H^1_{S_E}$,
the closed subscheme of $\prod_{a \in E} \P( S_a)$ whose subfunctor is
characterized by condition \eqref{quadraticEquations}, is cut out by
equations \eqref{BerndsOldBinomials}.
\end{proof}

In view of our general theory, Propositions \ref{p:graver} and
\ref{p:quad-binoms} show that Peeva and Stillman's determinantal
equations in \cite{PS} define the same scheme structure as the
binomial quadrics in \cite[Equation~(5.3)]{S}.  This question had been
left open in \cite{PS}.

\smallskip

It is instructive to examine Theorem \ref{c:determinants} in the case
of the toric Hilbert scheme $H_{S}^{1}$. Suppose the grading of $S$ is
positive, let $D \subseteq A_+$ be the set of Graver degrees and $E$
the very supportive set in \cite[Theorem 5.2]{PSlocal}. Then the toric
Hilbert scheme $H_{S}^{1}$ is defined by the quadratic binomials
\eqref{BerndsOldBinomials} on $D$ together with the Fitting equations
\eqref{e:determinantal}, where $e$ runs over $E$.  From this it
follows that the infinite sum in \cite[Definition 3.3]{PS} over all
degrees $e \in A_+$ can be replaced by the finite sum over $e \in E$.
The resulting finite set of determinantal equations still defines
$H_{S}^{1}$.

We now turn to the construction of the toric Chow morphism.  It was
conjectured in \cite[Problem 6.4]{S} that there
exists a natural morphism from the toric Hilbert scheme to a certain
inverse limit of toric GIT quotients, and this is what we shall now
construct.  In \cite{S} it was assumed that the action of $G$ on $\A
^{n}$ is the linearization of an action on projective space, or,
equivalently, that $(1,1,\ldots,1) \in M^\perp$, but this hypothesis is not
needed.  Our notation concerning toric varieties follows \cite{Cox}
and \cite{Ful}.  For compatibility with the standard toric variety
setting, one should take $k=\C $, although in fact the construction
below makes sense for any $k$.

In \eqref{e:homomorphismA->T} we identified $\Spec k[\x ,\x
^{-1}]$ with the torus $\T ^{n}$ of diagonal matrices acting on $\A
^{n}$.  Each Laurent monomial $x^{u}$ is thus a regular function on
$\T ^{n}$, and this identifies the lattice $\Z ^{n}$ in
\eqref{e:Apresentation} with the lattice of linear characters of $\T
^{n}$.  The sublattice $M$ consists of those characters which are
trivial on the subgroup $G\subseteq \T ^{n}$, so the grading group $A$
consists of linear characters of the group $G$.  In
particular, $S_{0}$ is the ring of $G$-invariants in $S$, so $\Spec
S_{0}$ is the affine quotient $\A ^{n}/G$.

The GIT quotient $\A ^{n}\slashsub{a} G$ in the sense of Mumford
\cite{Mum}, for a $G$-linearization of the trivial line bundle on $\A
^{n}$ using a $G$-character $a\in A$, is given as
\[ \A ^{n}\slashsub{a} G  \quad =  \quad 
\Proj \; \bigoplus_{r = 0}^\infty S_{rka}, \]
where $ka$ is multiple of $a$ for which the ring on the right-hand
side is generated in degrees $r=0,1$.  These GIT quotients,
including the affine quotient $\A ^{n}/G = \A ^{n}\slashsub{0}G$, are
toric varieties, whose description in terms of fans we
pause to review.  Let $N = \Hom (M,\Z)$ be the lattice dual to $M$, so
\eqref{e:Apresentation} yields an exact sequence
\begin{equation}\label{secondsequence}
0 \leftarrow \operatorname{Ext}^{1}(A,\Z ) \leftarrow N \leftarrow
\Z^{n} \leftarrow \Hom (A, \Z ) \leftarrow 0.
\end{equation} 
The map $N \leftarrow \Z^n$ supplies a tuple $\Pi = (v_{1},\ldots,v
_{n}) $ of distinguished vectors in $N$.  The lattice $\Z
^{n}$ in \eqref{secondsequence} is dual to the one in
\eqref{e:Apresentation}. The latter can be identified with
the set $\Z ^{\Pi }$ of functions $\Pi \rightarrow \Z$.  Then the
degree map $\,\Z ^{\Pi }\rightarrow A \,$ induces a map
\begin{equation}\label{e:map-phi}
\phi \, \colon\, \R ^{\Pi }\, \, \longrightarrow \,\, A\otimes _{\Z} \R .
\end{equation}
To each degree $a\in A$ is associated a regular subdivision $\Sigma
_{a}$, as in \cite{BGS}.  Namely, $\Sigma _{a}$ is the fan of cones in
$N$ spanned by subsets of the form $\sigma = f^{-1}(0) \subseteq \Pi
$, for functions $f\in (\R _{\geq 0})^{\Pi }\cap \phi ^{-1}(a)$.  The
fans $\Sigma =\Sigma _{a}$ arising this way are called {\it compatible
fans}.

With this notation, we have $X_{\Sigma _{a}} = \A ^{n}\slashsub{a} G$,
and we can identify $S$ with the Cox homogeneous coordinate ring
\cite{Cox} of $X_{\Sigma _{a}}$.  More precisely, for each $\sigma \in
\Sigma _{a}$ as above, let $x_{\sigma }$ be the product of the
variables $x_{i}$ with $v_{i}\in \Pi \setminus \sigma $.  Then the
semistable locus $U_{a} = \A ^{n}_{ss(a)}$ is the union of the
principal open affines $U_{x_{\sigma }}$, and the ring of
$G$-invariants $S[x_{\sigma }^{-1}]_{0}$ is the semigroup algebra of
$\sigma ^{\vee } \cap M$, so $U_{a}/G = X_{\Sigma _{a}}$.  It is
natural at this point to make the following definition.

\begin{definition}\label{d:integral}
A degree $a\in A_{+}$ is {\it integral} if the inclusion of convex
polyhedra 
\begin{equation} \label{IPvsLP}
\operatorname{conv} \bigl( \N^n \cap \deg^{-1}(a) \bigr) \,\,\subseteq \,\,
(\R_{\geq 0})^n \cap \phi^{-1}(a)
\end{equation}
is an equality, where $\phi $ is the map in \eqref{e:map-phi}.
\end{definition}

For every degree $a\in A_{+}$, the monomials $x_{\sigma }$ for $\sigma
\in \Sigma _{a}$ are just the square-free parts of all monomials
$x^{u}$ whose degree is a positive multiple of $a$.  Integrality of
$a$ means that every $x_{\sigma }$ already occurs as the square-free
part of a monomial of degree $a$.  Equivalently, $a$ is integral if
the semistable locus $U_{a}$ is equal to the complement of the closed
subset $V(J_{a})$, where $J_{a}$ is the ideal in $S$ generated by
$S_{a}$.  We remark that every degree $a\in A_{+}$ has some positive
multiple $k a$ which is integral, and that the fan $\Sigma _{a}$ is
the normal fan to the polyhedron on the right-hand side in
\eqref{IPvsLP}.

The set of all compatible fans, ordered by refinement, can be
identified with the poset of chambers in the Gale dual of $\Pi $, by
\cite[Theorem 2.4]{BGS}.  If $\Sigma'$ is a refinement of $\Sigma$
then the construction in \cite[\S 1.4]{Ful} gives a projective
morphism of toric varieties
\begin{equation}
\label{mapbetweentorics}
 X_{\Sigma'} \rightarrow  X_{\Sigma}. 
\end{equation}
The varieties $X_\Sigma$ for all compatible fans $\Sigma $ form an
inverse system.  Their inverse limit in the category of $k$-schemes is
called the {\em toric Chow quotient} and denoted by
\begin{equation}
\label{FirstInverseLimit} \A ^n \slashsub{C} G = \varprojlim \, \{
X_\Sigma : \text{$\Sigma $ a compatible fan in $N$} \}.
\end{equation}

\begin{example}
The fan $\Sigma _{0}$ giving the affine quotient $\A ^{n}/G$ is just
the cone $\R_{\geq 0}\Pi $.  If the grading is positive then $\Pi $
positively spans $N$ and all compatible fans are complete.  In this
case, the affine quotient $\A ^{n}/G = X_{\Sigma _{0}}$ is a point,
and the toric GIT quotients are projective.  At the opposite extreme,
if $A$ is finite, then $\Pi $ is a basis of $N \otimes \R$, the only
toric GIT quotient is the affine one, and our Chow morphism below
coincides with Nakamura's Chow morphism $\Hilb ^{G}(\A
^{n})\rightarrow \A ^{n}/G$.  \qed
\end{example}

The following theorem provides the solution to Problem 6.4 in \cite{S}.

\begin{theorem} \label{th:Prob64}
There is a canonical morphism
\begin{equation}\label{e:integralChow}
H^{1}_{S_{\integ (A)}} \rightarrow \A ^{n} \slashsub{C} G
\end{equation}
from the toric Hilbert scheme restricted to the set of integral
degrees to the toric Chow quotient $\A ^n \slashsub{C} G$, which induces an
isomorphism of the underlying reduced schemes.  In particular,
composing \eqref{e:integralChow} with the degree restriction morphism,
we obtain a canonical {\em Chow morphism} from the toric Hilbert
scheme to the toric Chow quotient
\begin{equation}
\label{toricChow} H^{1}_{S} \rightarrow \A ^{n} \slashsub{C} G.
\end{equation}
\end{theorem}

For the proof of Theorem \ref{th:Prob64} we need to recall some facts
about the $\Proj $ of a graded ring.  Let $T$ be an $\N$-graded
$k$-algebra, generated over $T_{0}$ by finitely many elements of
$T_{1}$, so $T = T_{0}[{\bf y}]/I$ for generators $y_{0},\ldots$,
$y_{m}$ of $T_{1}$ and a homogeneous ideal $I\subseteq T_{0}[{\bf
y}]$.  Then $\Proj T$ is the closed subscheme $V(I)$ of $\P ^{m}\times
\Spec T_{0}$.  Its functor $\fun{\Proj T}$ is defined as follows:
$\fun{\Proj T}(R)$ is the set of
homogeneous ideals $L\subseteq R\otimes T$ such that $L_{d}$ is a
locally free $R$-module of rank $1$ for all $d$.  In symbols,
$\Proj T = H^{1}_{T}$.  Setting $T^{(d)} = \oplus_{r} T_{rd}$, the
degree restriction morphism $\Proj T = H^{1}_{T}\rightarrow
H^{1}_{T^{(d)}} = \Proj T^{(d)}$ is an isomorphism.  More generally,
suppose that $T_{1}$ does not necessarily generate $T$, but that the
following weaker conditions hold:
\begin{itemize}
\item [(i)] $T$ is finite over the $T_{0}$-subalgebra $T'$ generated
by $T_{1}$; or equivalently,
\item [(ii)] there exists $d_{0}$ such that $T_{d+1} = T_{1}T_{d}$ for
all $d\geq d_{0}$.
\end{itemize}
In this case it remains true that the degree restriction morphism
$H^{1}_{T}\rightarrow \Proj T^{(d)}$ is an isomorphism for $d\geq
d_{0}$.  There is a canonical morphism $\Proj T^{(d)}\rightarrow \Proj
T'$, which is finite, but need not be an isomorphism.

One sees easily that $a\in A_{+}$ is integral if and only if the ring
\begin{equation}
S^{(a)} \underset{\text{def}}{=}\; \bigoplus _{r=0}^{\infty }S_{ra}
\end{equation}
satisfies conditions (i) and (ii) above.  In particular,  for
$a\in \integ (A)$ we have
\[
H^{1}_{S^{(a)}} = \A ^{n}\slashsub{a}G = X_{\Sigma _{a}},
\]
so restriction of degrees yields a canonical morphism $H^{1}_{S}
\rightarrow X_{\Sigma _{a}}$.

\begin{lemma}\label{l:a-vs-sigma}
The morphism $H^{1}_{S}\rightarrow X_{\Sigma }$ above depends only on
the compatible fan $\Sigma $ and not on the choice of degree $a$ with
$\Sigma _{a} = \Sigma $.  Moreover these morphisms commute with the
morphisms $X_{\Sigma '}\rightarrow X_{\Sigma }$ given by refinement of
fans as in (\ref{mapbetweentorics}).
\end{lemma}

\begin{proof}
We will describe the morphisms  $H^{1}_{S}\rightarrow X_{\Sigma }$
geometrically.  Because
$H^{1}_{S}$ represents the Hilbert functor, it comes with a universal
family $F\subseteq H^{1}_{S}\times _{k}\A ^{n}$ (where $\A ^{n} =
\Spec S$), and the group $G = \Spec k[A]$ acts on $F$ so
that the projections
\[
H^{1}_{S}\leftarrow F \rightarrow \A ^{n}
\]
are equivariant.  To the character $a$ of $G$ corresponds a GIT
quotient $F\slashsub{a}G = F_{ss}/G$.  The $a$-semistable locus
$F_{ss}$ is the preimage of $\A ^{n}_{ss}$, and we have induced
morphisms
\begin{equation}\label{e:geometric-morph}
H^{1}_{S} \leftarrow   F\slashsub{a}G  =
 F_{ss}/G \rightarrow \A^{n}_{ss}/G  =  X_{\Sigma }.
\end{equation}
Now, $F\slashsub{a}G$, considered as a scheme over $H^{1}_{S}$, is
just ${\Proj \bigoplus _{r=0}^{\infty } (\SS /\LL )_{ra}}$,
where $\SS $ is the sheaf of $A$-graded algebras ${\OO
_{H^{1}_{S}}\otimes _{k} S}$, and $\LL $ is the universal ideal
sheaf.  But $(\SS /\LL )_{ra}$ is locally free of rank $1$
over $\OO _{H^{1}_{S}}$ for all $r$, which implies
$F\slashsub{a}G\cong H^{1}_{S}$.  Hence in \eqref{e:geometric-morph}
there is a composite morphism $H^{1}_{S}\rightarrow X_{\Sigma }$,
which is easily seen to coincide with the degree restriction morphism.
The morphism in \eqref{e:geometric-morph} depends on $a$ only through
$\A ^{n}_{ss}$, which in turn depends only on $\Sigma _{a}$.
Furthermore, if $\Sigma '$ refines $\Sigma $, with corresponding
semi-stable loci $\A ^{n}_{ss'}$ and $\A ^{n}_{ss}$, then $\A
^{n}_{ss'}\subseteq \A ^{n}_{ss}$ and the morphism $X_{\Sigma
'}\rightarrow X_{\Sigma }$ is just the morphism $\A
^{n}_{ss'}/G\rightarrow \A ^{n}_{ss}/G$ induced by the inclusion.
This makes the lemma obvious.
\end{proof}

\begin{proof}[Proof of Theorem \ref{th:Prob64}]
We already saw  that restriction of degrees yields morphisms
\[
H^{1}_{S}\rightarrow H^{1}_{S_{\integ (A)}}\rightarrow \A ^{n}\slashsub{C} G.
\]
For every $R$, the natural map
\begin{equation}\label{e:int-nat}
H^{1}_{S_{\integ (A)}}(R)\rightarrow (\fun{\A ^{n}\slashsub{C} G})(R)
\end{equation}
is injective, since to give the restriction of $L\in H^{1}_{S}(R)$ to
integral degrees is the same as to give its image in
$H^{1}_{S^{(a)}}(R) = \fun{X_{\Sigma _{a}}}(R)$ for each integral $a$.
In general, a morphism of schemes $X\rightarrow Y$ induces an
isomorphism of reduced schemes $X_{\red }\rightarrow Y_{\red }$ if and
only if the natural map $\fun{X}(R)\rightarrow \fun{Y}(R)$ is an
isomorphism for all reduced rings $R$.  Hence it remains to show that
the map  \eqref{e:int-nat} is surjective when $R$ is reduced.

Suppose for each $a\in \integ (A)$ we are given $L_{a}\subseteq
{R\otimes S_{a}}$ with ${(R\otimes S_{a})}/L_{a}$ locally free of rank
1.  We may assume the $L_{a}$ are consistent in the following sense:
first, $\bigoplus _{r}L_{ra}$ is an ideal in $R\otimes S^{(a)}$ for
each $a$, so it represents a point of $\fun{X_{\Sigma _{a}}}(R)$, and
second, these points are compatible with the morphisms $X_{\Sigma
_{a}}\rightarrow X_{\Sigma _{b}}$ whenever $\Sigma _{a}$ refines
$\Sigma _{b}$.  Then we are to show:
\begin{equation}\label{e:consistent}
\text{for all $a,b\in \integ (A)$ and all $x^{u}$ with $\deg (u)=b-a$:
$x^{u}L_{a}\subseteq L_{b}$},
\end{equation}
so $L$ represents a point of $H^{1}_{S_{\integ (A)}}(R)$.

Let $D_{\Sigma } =\{a\in \integ (A):\Sigma _{a} = \Sigma \}$.  There
is a subdivision of $A\otimes _{\Z}\R$ into rational convex polyhedral
cones $C_{\Sigma }$ such that $D_{\Sigma }$ is the preimage of the
relative interior of $C_{\Sigma }$ via the canonical map $\psi \colon \integ
(A)\rightarrow A\otimes _{\Z}\R$.  The fan $\Sigma '$ refines $\Sigma
$ if and only if $C_{\Sigma }$ is a face of $C_{\Sigma '}$.  For a
given $\Sigma $, the ring $S_{\psi ^{-1}(C_{\Sigma })}$ can be
identified with the multi-homogeneous coordinate ring of $X_{\Sigma }$
with respect to the various line bundles $\OO (a)$ pulled back via
refinement homomorphisms $X_{\Sigma }\rightarrow X_{\Sigma _{a}}$.
Our consistency hypotheses amount to saying that \eqref{e:consistent}
holds whenever $\psi (a)$ and $\psi (b)$ both lie in a common cone
$C_{\Sigma }$.

Consider now the general case of \eqref{e:consistent}.  For any $d>0$,
the points $\psi (da+k\deg (u))$, for $0\leq k\leq d$, lie along the
line segment $l$ from $\psi (da)$ to $\psi (db)$.  For a suitably
chosen large $d$, every cone $C_{\Sigma }$ that meets $l$ will contain
at least one point $\psi (da+k_{i}\deg (u))$ with $da+k_{i}\deg (u)$
an integral degree.  Then \eqref{e:consistent} holds for each
consecutive pair of integral degrees $da+k_{i}\deg (u)$,
$da+k_{i+1}\deg (u)$ in this arithmetic progression.  Hence it holds
for $da$, $db$ and the monomial $x^{du}$, so $f\in L_{a}$ implies
$x^{du}f^{d}\in L_{db}$.  But $R\otimes S^{(b)}/(\bigoplus
_{r}L_{rb})$ is an $\N $-graded $R$-algebra, locally (on $\Spec R$)
isomorphic to a polynomial ring in one variable over $R$.  Hence it is
a reduced ring, that is, $\bigoplus _{r}L_{rb}$ is a radical ideal in
$R\otimes S^{(b)}$, and therefore $x^{u}f\in L_{b}$.
\end{proof}

Because the natural map in \eqref{e:int-nat} is always injective, our
proof of Theorem \ref{th:Prob64} gives a bit more.  Namely, if the
toric Chow quotient $\A ^{n}\slashsub{C} G$ happens to be reduced, then its
isomorphism with $(H^{1}_{S})_{\red }\subseteq H^{1}_{S}$ provides a
right inverse to the map in \eqref{e:int-nat}, showing that the latter
map is bijective.  Hence we have the following improvement.

\begin{corollary}\label{c:Chow-red}
If the toric Chow quotient
 $\A ^{n}\slashsub{C} G$ is reduced, then the morphism
$\,H^{1}_{S_{\integ (A)}}$ $ \rightarrow \A ^{n} \slashsub{C} G$ in Theorem
\ref{th:Prob64} is an isomorphism.
\end{corollary}

The toric Chow morphism is generally neither injective nor surjective,
see e.g.~\cite[Theorem 10.13]{StuBook}.  However, there is an
important special case, namely, the supernormal case, when it is
bijective, and in fact induces an isomorphism of the underlying
reduced schemes.  A degree $a \in A$ is called {\em prime} if there is
no variable $x_i$ which divides every monomial of degree $a$.  If
every degree is integral, ${\rm int}(A) = A$, then the sublattice $M$
of $\Z^n$ is said to be {\em unimodular}.  If every
prime degree is integral, the $A$-grading is said to be {\em
supernormal}.  This terminology is consistent with \cite{HMS}.

\begin{corollary} \label{SuperNormalIso} 
If the $A$-grading of $S$ is supernormal then the
toric Chow morphism \eqref{toricChow} induces an isomorphism of
reduced schemes $(H^{1}_{S})_{\red }\rightarrow (\A ^{n}\slashsub{C}
G)_{\red }$.  If in addition, $\A ^{n}\slashsub{C} G$ is reduced, then
the toric Chow morphism is an isomorphism.
\end{corollary}

\begin{proof}
We need only show that if $D$ is the set of prime degrees, then the
degree restriction morphism $H^{1}_{S}\rightarrow H^{1}_{S_{D}}$ is an
isomorphism.  Fix a degree $a \in A \backslash D$ and let $x^u$ be the
greatest common divisor of all monomials of degree $a$.  Then $a' = a
- \deg(u)$ is a prime degree, and multiplication by $x^u$ defines an
$R$-module isomorphism between $R \otimes S_{a'}$ and $R \otimes
S_{a}$, for every $k$-algebra $R$.  For any element $I \in
H^{1}_S(R)$, we have $I_a = x^u \cdot I_{a'} $ and hence the
restriction map $ H^{1}_S(R) \rightarrow H^{1}_{S_D}(R) $ is
injective. But it is also surjective because every element $J$ in
$H^{1}_{S_D}(R) $ lifts to an element of $H^{1}_{S}(R) $ by setting
$J_a = x^u \cdot J_{a'} $ for degrees $a \in A \backslash D$.
\end{proof}

\begin{example}
Give $S = \C[x_1,x_2,x_3,x_4]$ the $\Z^2$-grading $ \deg(x_1) =
\deg(x_2) = (1,0)$, $\deg(x_3) = (0,1)$, $\deg(x_4) = (2,1)$.  The
configuration $\Pi \subseteq N$ can be represented by the four vectors
$\{(-1,1),\, (1,1),\, (0,1),\, (0,-1) \}$ in $\Z ^{2}$.  There is a
unique finest compatible fan $\Sigma $, and $\A ^{2}\slashsub{C} G =
X_{\Sigma }$ is a smooth projective toric surface.  The prime degrees
are $(\alpha ,\beta )\in \N ^{2}$ with $\alpha \geq 2\beta $.  The
integral degrees are those for which $\alpha \geq 2\beta $ or $\alpha
$ is even.  Hence this example is supernormal, but not unimodular.  By
Corollary \ref{SuperNormalIso}, its toric Hilbert scheme is isomorphic
to $X_{\Sigma }$.  \qed
\end{example} 

We remark that Corollary \ref{SuperNormalIso} can be used to give an
alternative and more conceptual proof of Theorem 1.2 in \cite{HMS},
which states that in the supernormal case, the $\T ^{n}$-fixed points
on the toric Hilbert scheme (``virtual initial ideals'') are in
natural bijection with the $\T ^{n}$-fixed points on the toric Chow quotient
(``virtual chambers'').

We next discuss the distinguished component of the toric Hilbert
scheme $H^{1}_S$.  Here we will fix $k=\C $ and describe the
distinguished component in set-theoretic terms.  The distinguished
point $I_{M}\in H^{1}_{S}(\C )$ is the ideal of the closure of $G$.
More generally, the ideal $I$ of the closure in $\A ^{n}$ of any
$G$-coset $G\cdot \tau \subseteq \T ^{n}$ is a point of $H^{1}_{S}(\C
)$.  In fact, $\T ^{n}$ acts on $H^{1}_{S}$ and the $\T ^{n}$-orbit of
$I_{M}$ consists of all such ideals $I$.  Moreover, $I$ is the only
point of $H^{1}_{S}(\C )$ for which $V(I)$ contains $G\cdot \tau $.
Now $\T ^{n}/G$ is the open torus orbit in each of the toric varieties
$X_{\Sigma }$, and so is naturally embedded as an open set $U$ in the
inverse limit $\A ^{n}\slashsub{C}G = \varprojlim_{\Sigma } X_{\Sigma
}$.  The observations above show that the toric Chow morphism
restricts to a bijection from the preimage of $U$ in $H^{1}_{S}$ to
$U$.  Hence the preimage of $U$ is an irreducible open subset of
$H^{1}_{S}$, and its closure is an irreducible component of
$H^{1}_{S}$, which we call the {\it coherent component}.

The closure of $U$ in  $\,\A ^{n}\slashsub{C}G \,$ is the
toric variety $X_{\Delta }$ defined by the common refinement $\Delta $
of all  compatible fans $\Sigma _{a}$.  Thus we have a canonical
morphism from the coherent component of $H^{1}_{S}$ to $X_{\Delta }$.
For a supernormal grading, it is an isomorphism.

\begin{example}
There is a nice connection between toric Hilbert schemes and recent
work by Brion \cite{Bri}. Consider the
(Grothendieck)  Hilbert
scheme associated with the diagonal embedding $X \rightarrow X \times
X$ of a projective variety $X$.
Brion shows that, if $X$ is a homogeneous
space, then the diagonal is a smooth point on a unique component of
the Hilbert scheme, the associated Chow morphism is an isomorphism,
and all degenerations of the diagonal in $X \times X$ are reduced and
Cohen-Macaulay.
Our theory implies corresponding results when $X$ is a unimodular
toric variety \cite[\S 6]{BPS}.  Indeed, by Proposition 6.1 in
\cite{BPS}, the Hilbert scheme of $X$ in $X \times X$ is the toric
Hilbert scheme for $S = \C[x_1,\ldots, x_n,y_1,\ldots,y_n]$ graded via
a unimodular Lawrence lattice $M \subseteq \Z^{2n}$.  The Lawrence
ideal $I_M$ is the distinguished point $X$; the unique component it
lies on is the distinguished component defined above, and its toric
degenerations are reduced and Cohen-Macaulay by \cite[Thm.~1.2
(b)]{BPS}. \qed
\end{example}

Variants of our results on $H^{1}_{S}$ apply to the {\em $m$-th toric
Hilbert scheme} $H^{m}_{S}$ defined by the Hilbert function
\[
h(a) = \min (m, \rank _{k}(S_{a})).
\]
The $m$-th toric Hilbert scheme is a common generalization of two
objects of recent interest in combinatorial algebraic geometry,
namely, the Hilbert scheme of $m$ points in affine $n$-space (the case
$A = \{0\}$) and the toric Hilbert scheme (the case $m=1$).  Again,
$H^{m}_S$ has a distinguished {\em coherent component}, and it admits
a natural morphism to a certain Chow variety.

We briefly outline how to extend our constructions to the case $m>1$.
The appropriate Chow quotient is the inverse limit of symmetric powers
$\varprojlim_{\Sigma }\Sym ^{m} X_{\Sigma }$.  There is a Chow morphism,
which factors as
\begin{equation}\label{e:chow-factor}
H^{m}_{S}\rightarrow \prod _{\Sigma }\Hilb ^{m}( X_{\Sigma
})\rightarrow \varprojlim_{\Sigma }\Sym ^{m} X_{\Sigma }.
\end{equation}
Here $\Hilb ^{m}(X_{\Sigma })$ is the Hilbert scheme of $m$ points in
$X_{\Sigma }$. For
sufficiently general $a$ such that $\Sigma =\Sigma _{a}$, we have
$\Hilb ^{m}(X_{\Sigma }) = H^{m}_{S^{(a)}}$, and the Chow morphism is
degree restriction composed with the usual Chow morphisms $\Hilb
^{m}(X_{\Sigma })\rightarrow \Sym ^{m} X_{\Sigma }$.  The analog of
Theorem \ref{th:Prob64} no longer holds, however.

The coherent component is the unique component of $H^{m}_{S}$ which
maps birationally on the open subset $\Sym ^{m}(\T ^{n}/G)$ of the
Chow quotient.  A typical point in the coherent component is the ideal
of the union of $m$ closures in $\A ^{n}$ of $G$-cosets in $\T ^{n}$.
One difference from the $m=1$ case is that for $m>1$, not every ideal
of a union of $G$-coset closures belongs to $H^{m}_{S}$.  If the $m$
cosets are specially chosen, it can happen that the monomials of some
degree $a$ define fewer than $h(a)$ linearly independent functions on
them, a possibility that does not occur when $m=1$.

\section{Constructing other Hilbert schemes}
\label{sec6}

The theory of graded $k$-modules with operators developed in
Section 2 allows us to construct many interesting Hilbert schemes in
addition to the multigraded Hilbert schemes of Theorem \ref{t:multi-Hilb}. 
This final section lists some noteworthy examples.

\subsection{Partial multigraded Hilbert schemes}

Take $S$ an $A$-graded polynomial ring, as before, $T=S_{D}$ for any
subset $D \subset A$, and any function $h : D \mapsto \N$.  The
operators are multiplications by monomials.  Using Remark
\ref{r:Fclosed}, we see that the ``monomial ideals'' in the system
$(T,F)$ are just the restrictions to $S_{D}$ of monomial ideals in
$S$.  This given, it is easy to see that the analog of Proposition
\ref{p:Diane-prime} holds, and the proof of Theorem \ref{t:multi-Hilb}
goes through to show that $H^h_{T} = H^{h}_{S_{D}}$ is represented by
a quasiprojective scheme, projective if the grading is positive.  (We
actually used this result already in Section \ref{sec5} when we
implicitly assumed that the integral degrees Hilbert functor
$H^{1}_{S_{\integ (A)}}$ is represented by a scheme.)  Whenever $D
\subset E \subset A$, we have a degree restriction morphism
$H^h_{S_D} \rightarrow H^h_{S_E}$.  Such Hilbert schemes occur
naturally in parametrizing subschemes of a toric variety $X$.  Here
$S$ is the homogeneous coordinate ring \cite{Cox} of $X$, the grading
group $A$ is the Picard group of $X$ and $D$ is a suitable translate
of the semigroup of ample divisors in $A$.  This application is
currently being studied by Maclagan and Smith \cite{MS}.

\subsection{Quot schemes}

Take $T = S^r/M$ to be a finitely generated $A$-graded module over the
$A$-graded polynomial ring $S$.  The scheme representing $H^{h}_{T}$
is a {\it Quot scheme}, used in algebraic geometry to parametrize
vector bundles or sheaves on a given scheme.  The arguments in Section
\ref{sec3} extend easily to show that the multigraded Quot scheme is
always a quasiprojective scheme over $k$.  The special case $A = 0$,
$M =0 $ is already interesting: here $H_{T}^{h}$ parametrizes artinian
$S$-modules with $r$ generators having length $m=h(0)$.  For $n=2$,
this scheme is closely related to the space ${\mathcal M}(m,r)$ of
Nakajima, which plays the lead character in his work on Hilbert
schemes of points on surfaces; see \cite[Chapter 2]{Nak}.

\subsection{Universal enveloping algebras }

Let $\gg = \bigoplus _{a}\gg_{a}$ be an $A$-graded Lie algebra over
$k$, free and finitely generated as a $k$-module.  Take $T = U(\gg)$
to be its universal enveloping algebra, or, more generally, if $c\in
\gg$ is a central element of degree $0$, take $T$ to be the reduced
enveloping algebra $U_{c}(\gg) = U(\gg)/I$, where $I$ is the two-sided
ideal $\langle c-1 \rangle$.  Then $T$ is an $A$-graded associative
algebra.  Taking our system of operators $F$ to be generated by the
left multiplications by elements of $\gg$, or the right
multiplications, or both, we obtain Hilbert functors $H^{h}_{T}$
parametrizing homogeneous left, right or two-sided ideals with Hilbert
function $h$.  Similarly, there are Quot functors $H^{h}_{T}$, when
$T$ is a finitely generated module over $U(\gg)$ or $U_{c}(\gg)$.

All these functors are represented by quasiprojective schemes over
$k$.  To see this, recall the {\it arithmetic filtration} of $U =
U(\gg)$ or $U_{c}(\gg)$ given by $\alpha ^{i}U = \sum _{j\leq
i}\gg^{j}$.  The associated graded algebra $\gr _{\alpha }U = S$ is an
$A$-graded commutative polynomial ring in variables
$x_{1},\ldots,x_{n}$ forming a homogeneous $k$-basis of $\gg$ (or of
$\gg/kc$).  If $I\subseteq U$ is a left, right or two-sided ideal,
then $\gr _{\alpha }I$ is an ideal in $S$.  When the ground ring is a
field, we may define the {\it initial ideal} $\init (I) = \init (\gr
_{\alpha }I)$.  Fixing a Poincar\'e-Birkhoff-Witt basis in $U$, the
basis elements corresponding to standard monomials for $\init (I)$ form
a basis of $U/I$.  Given these observations, it is easy to adapt the
proof of Theorem \ref{t:multi-Hilb} in Section \ref{sec3} to show that
in this more general setting, $H^{h}_{T}$ is still a closed subscheme
of the quasiprojective scheme $H^{h}_{T_{D}}$, where $D$ is a finite
supportive set for $S$ and $h$.  It can also be shown, using Gr\"obner
basis theory for $U$, that $H^{h}_{T}$ is isomorphic to
$H^{h}_{T_{D}}$ if $D$ is very supportive for $S$ and $h$.

An interesting example is the {\it Weyl algebra} $W = k \langle
x_1,\ldots,x_n, \partial_1,\ldots,\partial_n \rangle $, which is the
reduced enveloping algebra of a Heisenberg algebra $\gg$.  Any
$A$-grading of $k[\x ]$ extends to a grading of $W$ with $\deg
(\partial_i) = -\deg (x_{i}) $.  Then we have a Hilbert scheme
$H^{h}_{W}$ parametrizing homogeneous left ideals with Hilbert
function $h$ in the Weyl algebra.  It would be interesting to relate
these Hilbert schemes to the work of Berest and Wilson \cite{BW} in
the case $n = 1$.  We note that, unlike in the case of the polynomial
ring $S = k[\x ]$, the finest possible grading, $A = \Z^n$, $ \deg
(x_i) = - \deg(\partial_i) = e_i$, gives rise to highly non-trivial
Hilbert schemes $H^{h}_{W}$.  Namely, $H^{h}_{W}$ parametrizes all
{\it Frobenius ideals} with Hilbert function $h$.  Frobenius ideals
appear in Gr\"obner-based algorithms for solving systems of linear
partial differential equations \cite[\S 2.3]{SST}.

\subsection{Quivers, posets and path algebras}

Here is a nice example where the set $A$ of ``degrees'' is not a group.
Fix a finite poset $Q$.  We identify $Q$ with its Hasse diagram and we
interpret it as an acyclic quiver.  Let $ T = k Q$ denote its {\it
path algebra}.  This is the free $k$-module spanned by all directed
paths in $Q$ modulo the obvious concatenation relations.  We take $A$
to be the set of all intervals $[u,v]$ in the poset $Q$.  Then the
``graded component'' $T_{[u,v]}$ is the $k$-vector space with basis
consisting of all chains from $u$ to $v$ in the poset $Q$.  Fixing a
Hilbert function $h : A \rightarrow \N$, we get the scheme $H^h_T$
which parametrizes all homogeneous quotients of the path algebra $k Q$
modulo a certain number of linearly independent relations for each
interval $[u,v]$.  The case when $h$ attains only the values $0$ and
$1$ deserves special attention.  In this case the scheme $H^h_T$ is
binomial, and it parametrizes the {\it Schurian algebras} in the sense
of \cite{CPPRT}.  When $h = 1$ is the constant one function then we
get a non-commutative analogue to the toric Hilbert scheme.  The
distinguished point on this Hilbert scheme $H^1_T$ is the {\it
incidence algebra} of the poset $Q$, and it would be interesting to
study its deformations from this point of view.  Question: What is the
smallest poset $Q$ for which the scheme $H^1_T$ has more than one
component?

\subsection{Exterior Hilbert schemes }

The use of the exterior algebra as a tool for (computational)
algebraic geometry has received considerable attention in recent
years; see e.g.~\cite{EFS}.  Let $T$ denote the exterior algebra in
$n$ variables $x_1,\ldots,x_n$ over our base ring $k$, again with a
grading by an abelian group $A$.  As in the preceding example, $T$ is
of finite rank over $k$, so the Hilbert scheme $H^{h}_{T}$ is a closed
subscheme of a product of Grassmann schemes (as are the associated
Quot schemes).  The torus fixed points on such Hilbert schemes are
precisely the simplicial complexes on $\{x_1,\ldots,x_n\}$, and it
would be interesting to study the combinatorial notion of ``shifting
of simplicial complexes'' in the framework of exterior Hilbert
schemes.

\end{document}